\documentclass[a4paper,12pt]{article}

\usepackage{amsmath,amssymb,amsfonts,graphicx}

\usepackage{epsfig}
\usepackage{color}
\usepackage{comment}
\usepackage{ulem}
\usepackage{mathdots}
\usepackage{graphics}
\usepackage{epstopdf}
\usepackage{mathdots}
\usepackage{epsfig}
\usepackage{verbatim}
\usepackage{color}
\usepackage{array}
\usepackage{url}
\usepackage{textcomp}

\newcolumntype{P}[1]{>{\centering\arraybackslash}p{#1}}

\def\sq{ { {\color{blue} \hfill\rule{1.5mm}{1.5mm} } } }

\newtheorem{theorem}{Theorem}[section]
\newtheorem{lemma}{Lemma}[section]
\newtheorem{definition}{Definition}[section]
\newtheorem{proposition}{Proposition}[section]
\newtheorem{corollary}{Corollary}[section]
\newtheorem{remark}{Remark}[section]

\newtheorem{example}{Example}[section]
\newtheorem{assumption}{Assumption}[section]
\newtheorem{procedure}{Procedure}[section]

\newcommand{\bA}{{\bf A}}
\newcommand{\bB}{{\bf B}}
\newcommand{\bC}{{\bf C}}
\newcommand{\bD}{{\bf D}}
\newcommand{\bE}{{\bf E}}

\newcommand{\bS}{{\bf S}}

\newcommand{\bK}{{\bf K}}

\newcommand{\bM}{{\bf M}}
\newcommand{\bN}{{\bf N}}
\newcommand{\bI}{{\bf I}}
\newcommand{\bH}{{\bf H}}
\newcommand{\bP}{{\bf P}}

\newcommand{\bQ}{{\bf Q}}

\newcommand{\bR}{{\bf R}}

\newcommand{\bZ}{{\bf Z}}
\newcommand{\bX}{{\bf X}}
\newcommand{\bx}{{\bf x}}
\newcommand{\by}{{\bf y}}
\newcommand{\bu}{{\bf u}}
\newcommand{\bV}{{\bf V}}
\newcommand{\bU}{{\bf U}}

\newcommand{\bfe}{{\bf e}}

\newcommand{\bn}{{\bf n}}

\newcommand{\bg}{{\bf g}}

\newcommand{\bh}{{\bf h}}

\newcommand{\bv}{{\bf v}}

\newcommand{\bz}{{\bf z}}

\newcommand{\bff}{{\bf f}}
\newcommand{\bfz}{{\mathbf 0}}

\newcommand{\cP}{ {\cal P} }
\newcommand{\cQ}{ {\cal Q} }

\newcommand{\Si}{ \boldsymbol{\Sigma} }

\newcommand{\bPhi}{ \boldsymbol{\Phi} }

\newcommand{\bLa}{\boldsymbol{\Lambda}}

\title{
   Balanced truncation for linear switched systems
}
\author{\\
I.V. Gosea \footnotemark[2]~ ,~  M. Petreczky \footnotemark[3]~,
   A.C. Antoulas \footnotemark[4]  \ \footnotemark[2] ,
      C. Fiter \footnotemark[5]
}

\begin{document}

\maketitle

\renewcommand{\thefootnote}{\fnsymbol{footnote}}
\footnotetext[2]{
Data-Driven System Reduction and Identification Group, Max Planck Institute for Dynamics of Complex Technical Systems, Sandtorstrasse 1, 39106, Magdeburg, Germany
({\tt gosea@mpi-magdeburg.mpg.de})
}

\footnotetext[3]{Centre de Recherche en Informatique, Signal et Automatique de Lille (CRIStAL),  UMR CNRS 9189, CNRS, Ecole Centrale de Lille, France
({\tt mihaly.petreczky@ec-lille.fr})
}

\footnotetext[4]{
   Department of Electrical and Computer Engineering, 
   Rice University, 6100 Main St, MS-366, Houston, TX 77005, USA 
   ({\tt aca@rice.edu})
}

\footnotetext[5]{
    CNRS CRIStAL UMR 9189, Universit\'e de Lille 1, Sciences et Technologies, 59651 Villeneuve d\textquotesingle Ascq, France  
    ({ \tt christophe.fiter@univ-lille1.fr})}

\renewcommand{\thefootnote}{\arabic{footnote}}

\begin{abstract} 

We propose a model order reduction approach for balanced truncation of linear switched systems. Such systems switch among a finite number of linear subsystems or modes.

We compute pairs of controllability and observability Gramians  corresponding to each active discrete mode by solving systems of coupled Lyapunov equations.
Depending on the type, each such Gramian corresponds to the energy associated to all possible switching scenarios that start or, respectively end, in a particular operational mode.

In order to guarantee that hard to control and hard to observe states are simultaneously eliminated, we construct a transformed system, whose Gramians are equal and diagonal. Then, by truncation, directly construct reduced order models. One can show that these models preserve some  properties of the original model, such as stability and that it is possible to obtain error bounds relating the observed output, the control input and the entries of the diagonal Gramians.

\end{abstract}

\section{Introduction}

In recent years, the need for accurate mathematical modeling of physical and artificial processes for simulation and control has been steadily increasing. To cope with it, inclusion of more detail at the modeling stage is required, which inevitably leads to analyzing  larger-scale, more complex dynamical systems. Such high dimensional systems are often linked  to spatial discretization of underlying time-dependent coupled partial differential equations (PDE).

In broad terms, model order reduction (MOR) is concerned with finding efficient computational prototyping tools to replace such complex and large models by simpler and smaller models that capture their dominant characteristics. Such reduced order models (ROM) could be used as efficient surrogates for the original model, replacing it as a component in larger simulations. For details on different MOR techniques, we refer the readers to the book \cite{bookAA} and to the surveys \cite{baur14,benner15}.

Hybrid systems are a class of nonlinear systems which result from the interaction of continuous time dynamical subsystems with discrete events. These systems are hence described by both discrete and continuous states, inputs and outputs. 
The transitions between the discrete states may result in a jump in the continuous internal variable. The
discrete dynamics is determined by a finite-state deterministic automaton equipped with outputs (the so-called Moore automaton).

Switched systems constitute a subclass of hybrid systems, in the sense that the discrete dynamics is simplified, i.e. any discrete state transition is allowed and the set of discrete events coincides with the set of discrete states.

A switched system is a dynamical system that consists of
a finite number of subsystems and a logical rule that
orchestrates switching between these subsystems.
These subsystems  or discrete modes are usually described by a collection of  differential or difference equations. The discrete events interacting with the subsystems are governed by a piecewise continuous function, i.e. the switching signal.

One can classify switched systems based on the dynamics of
their subsystems, for example continuous-time or discrete-time, linear or nonlinear and so on. In this work we analyze continuous-time linear switched systems (LSS) with reset maps (or coupling/switching matrices). The latter term refers to matrices that scale the continuous state at the switching times.

Hybrid and switched systems represent useful models for distributed embedded systems design where discrete controls are routinely applied to continuous processes. In particular, switched systems have applications in control of mechanical  and aeronautical systems, power converters and also in the automotive industry. In some cases, the complexity of verifying and assessing general properties of these systems
is very high so that the use of these models is limited in
applications where the size of the state space is large. A useful tool for dealing with such complexity is MOR. For a detailed characterization of theses classes of dynamical systems, we refer the readers to the books \cite{book1LSS,book2LSS,book3LSS,book4LSS}.

In the past years, hybrid and switched systems have received increasing attention in the scientific community, which can be partly explained by the fast development of the switch control research area (see \cite{ndfg03,vcl07,llmk11}). In this context, adaptive control techniques based on switching between different controllers are used to achieve stability and improve transient response. 

The study of the properties of hybrid and switched systems includes such topics as stability (see \cite{dri02,vcl07,book2LSS}), realization including observability/controllability (see \cite{petre11,ptt13}), analysis of switched DAE's (see \cite{mehrmann09,tr12}) and numerical solutions (see \cite{mehrmann08}).
Considerable attention has been dedicated in recent years to the problem of MOR for linear switched systems. A class of methods that involves matching of generalized Markov parameters (known also as time domain Krylov methods) has been discussed in \cite{bastug16,thesisB}; $\mathcal{H}_{\infty}$ type of reduction methods were developed in \cite{zsbw08,bmb11,zcw14}. We mention some publications that are focused on the reduction of discrete LSS, such as \cite{zbs09} and \cite{bmb12}.

A very prolific MOR method that has been continuously developed over the years is balanced truncation (BT). It was initially introduced in the systems and control theory in \cite{Mo81,ps82}. The main idea behind BT is to transform a dynamical system to a balanced form defined in such a way that appropriately chosen controllability and observability
Gramians are equal and diagonal. Then a reduced-order model is computed by truncating the states corresponding to the small diagonal elements of the Gramians. For more details on BT especially from a practical point of view (i.e. application to large scale systems, solution of Lyapunov equations etc.), see \cite{ms05,bs17}.

For switched and hybrid systems, techniques that are based on balancing (or of Gramian based derivation of it) have been considered in the following: \cite{glw06,ch09,bgmb10,sw11,sw12,mtc12,pwl13,pp16}.

For LSS, it may happen that some state components are difficult to reach and observe in some of the modes yet easy to reach and observe in others. In that case, deciding how to truncate the state variables and obtain a meaningful ROM is not trivial. Under general conditions, there is no basis of the state space so that all modes of the LSS are in balanced form (having equal, diagonal Gramians). This problem was addressed in \cite{mtc12}. There, very restrictive necessary and sufficient conditions are derived for the existence of such basis.

We are interested in the situation for which a common transformation that simultaneously balances all the subsystems of the LSS is either not known or it does not exist. We construct a family of transformations, each for a specific mode, based on appropriately defined Gramians. Concerning stability, this concept extends to the existence of multiple Lyapunov functions (Chapter 3.1 in \cite{book1LSS}). Another key factor is that we consider a sufficient slow switching to impose certain properties such as stability, error and energy bounds. Hence the assumption of a minimum dwell time (the duration of time for which the system is active in a particular mode) which is chosen depending on the context or application.

The proposed method is centered around the definition of new type of Gramians for LSS that resemble the definitions previously encountered for the case of bilinear and stochastic systems (see \cite{zl02,bdc17}). Some of the results presented in this work are extended from \cite{pwl13}. There, the Gramians are defined as solutions of systems of linear matrix inequalities (LMI).

The paper is organized as follows; in the second section, we introduce continuous-time linear switched systems in a formal way. Furthermore, we provide a characterization of input-output mappings in time domain corresponding to such systems. Section 3 describes the procedure of constructing infinite energy Gramians for the simplified case with only two discrete modes. 
Next, in Section 4 we provide a system theoretic interpretation of such Gramians (for the general case with D modes). Furthermore, we formally introduce the balancing algorithm followed by the MOR step, i.e. the truncation. A measure of the quality of approximation by reduction is provided by means of a error bound. Additionally, we investigate the possibility of preserving system theoretic properties such as stability, for the reduced system. Section 5 is designated for the numerical experiments while a summary of the findings and the conclusion are presented in Section 6.

\section{Linear switched systems} \label{sec:switched}

\begin{definition}
	A continuous time linear switched system (LSS) is a control system of the form:
	\begin{equation}\label{LSS_def0}
	\Si: \begin{cases}
	\bE_{\sigma(t)} \dot{\bx}(t) = \bA_{\sigma(t)} \bx(t) + \bB_{\sigma(t)}\bu(t), \ \ \bx(t) = \bx_0, \\
	\by(t) = \bC_{\sigma(t)} \bx(t),
	\end{cases}
	\end{equation}
	where $\Omega = \{1,2,\ldots,D\}, \  D > 1$, is a set of discrete modes, $\sigma(t)$ is  the switching signal, $\bu$ is  the input, $\bx$ is the state, and $\by$ is the output.

	The system matrices  $\bE_q, \bA_q \in \mathbb{R}^{n_q \times n_q}, \ \bB_q \in \mathbb{R}^{n_q \times m}, \ \bC_q \in \mathbb{R}^{p \times n_q}$, where $q \in \Omega$, correspond to the linear system active in mode $q \in \Omega$, and $\bx_0$ is the initial state. We consider the $\bE_q$ matrices to be invertible. Furthermore, the transition from one mode to another is made via the so called switching or coupling matrices $\bK_{q_1,q_2} \in \mathbb{R}^{n_{q_2} \times n_{q_1}}$ where  $q_1,q_2 \in \Omega$. 
\end{definition}

\begin{remark}
	{\rm The case for which the coupling is made between identical modes is excluded, Hence, when $q_1 = q_2 = q$, consider that the coupling matrices are identity matrices, i.e. $\bK_{q,q} = \bI_{n_q}$.}
\end{remark}
\noindent
The notation $\Si = (n_1,n_2,\ldots,n_D, \{(\bE_q,\bA_q, \bB_q, \bC_q)| q \in \Omega\},\{ \bK_{q_{i},q_{i+1}} | q_i, q_{i+1} \in \Omega \}, \bx_0)$ is used as a short-hand representation for LSS's described by the equations in (\ref{LSS_def0}). The vector $\bn = \left( \begin{array}{cccc}
n_1 & n_2 & \cdots & n_D
\end{array} \right)$ is the dimension (order) of $\Si$. 

The restriction of the switching signal $\sigma(t)$ to a finite interval of time $[0, T]$ can be interpreted as finite sequence of elements of $\Omega \times \mathbb{R}_{+}$ of the form:
$$
\nu(\sigma) = (q_1,t_1) (q_2,t_2) \ldots (q_k,t_k),
$$
where $q_1, \ldots , q_k \in \Omega$ and $0<t_1 < t_2 < \cdots < t_k  \in \mathbb{R}_{+}, \ t_1 + \cdots + t_k = T$, such that for all $ t \in [0, T]$ we have:
$$
\sigma(t) = \left\{ \begin{array}{l} q_1 \ \  \text{if} \ \ t \in [0,t_1], \\
q_2 \ \  \text{if} \ \  t \in (t_1,t_1+t_2], \\ 
\ldots \\
q_i \ \  \text{if} \ \ t \in (t_1+\ldots+t_{i-1},t_1+\ldots+t_{i-1}+t_i], \ \text{for} \ 2 \leqslant i \leqslant k.
\end{array} \right.
$$
In short, by denoting $T_i := t_1+\ldots+t_{i-1}+t_i, \  T_0 :=0, \  T_k := T$, write
\begin{equation}\label{switch_signal}
\sigma(t) = \begin{cases}  q_1 \ \  \text{if} \ \ t \in [0,T_1], \\
q_i \ \  \text{if} \ \ t \in (T_{i-1},T_i], \ i > 2. \end{cases}
\end{equation}

The linear system which is active in the $q_i^{th}$ mode of $\Si$ is denoted with $\Si_{q_i}$ and it is described by (where $\bx_{q_i}(t) = \bx(t), \ t \in (T_{i-1},T_i]$)
\begin{equation}\label{LSS_def2}
\Si_k: \begin{cases}
\bE_{q_i} \dot{\bx}_{q_i}(t) = \bA_{q_i}\bx_{q_i}(t) + \bB_{q_i}\bu(t), \ \   \\
\by(t) = \bC_{q_i} \bx_{q_i}(t). \ \ 
\end{cases}
\end{equation}

Denote by $PC(\mathbb{R}_{+},\mathbb{R}^n)$, $P_c(\mathbb{R}_{+},\mathbb{R}^n)$, the set of all piecewise-continuous, and piecewise-constant functions, respectively.
\begin{definition}
	A tuple $(\bx,\bu,\sigma,\by)$, where $\bx:\mathbb{R}_{+} \rightarrow \bigcup_{i=1}^{D} \mathbb{R}^{n_i}$, $\bu \in PC(\mathcal{R}_{+},\mathbb{R}^{m}),$ $ \sigma \in P_c(\mathbb{R}_{+},\Omega), \by \in PC(\mathbb{R}_{+},\mathbb{R}^{p})$ is called a solution, if the following conditions simultaneously hold:
	
	\begin{enumerate}
		\item The restriction of $\bx(t)$  to $(T_{i-1},T_{i}]$ is differentiable, and satisfies $\bE_{q_i} \dot \bx (t) =\bA_{q_i}\bx(t)+\bB \bu(t)$.
		\item  Furthermore, when switching from mode $q_{i}$ to mode $q_{i+1}$ at time $T_i$, the following holds
		\begin{equation*}
		\bE_{q_{i+1}} \displaystyle \lim_{t \searrow T_i} \bx(t) = \bK_{q_i,q_{i+1}}  \bx(T_i).
		\end{equation*}  
		\item  Moreover, for all $t \in \mathbb{R}$, $\by(t)=\bC_{\sigma(t)} \bx(t)$ holds.
	\end{enumerate}
\end{definition}

The switching matrices $\bK_{q_{i},q_{i+1}}$ allow having different dimensions for the subsystems active in different modes. For instance, the pencil $(\bA_{q_i},\bE_{q_i}) \in \mathbb{R}^{n_{q_i} \times n_{q_i}}$, while the pencil $(\bA_{q_{i+1}},\bE_{q_{i+1}}) \in \mathbb{R}^{n_{q_{i+1}} \times n_{q_{i+1}}}$ where the values $n_{q_i}$ and $n_{q_{i+1}}$ need not be the same. If the $\bK_{q_{i},q_{i+1}}$ matrices are not explicitly given, it is considered that they are identity matrices.

The input-output behavior of an LSS system can be described in time domain using the mapping $\bff(\bu,\sigma)$. This particular map can be written in \textit{generalized kernel representation} (as suggested in \cite{pvs10}) using the unique family of analytic functions: $\bg_{q_1,\ldots,q_k}: \mathbb{R}_{+}^k \rightarrow \mathbb{R}^p$ and $\bh_{q_1,\ldots,q_k}: \mathbb{R}_{+}^k \rightarrow \mathbb{R}^{p \times m}$ with $q_1,\ldots,q_k \in \Omega, \ k \geqslant 1$ such that for all pairs $(\bu,\sigma)$ and for $T = t_1+t_2+\cdots+t_k$ we can write:
\small
\begin{align}
& \bff(\bu,\sigma)(t) = \bg_{q_1,\ldots,q_k}(t_1,...,t_k)+  \displaystyle \sum_{i=1}^{k} \int_{0}^{t_i} \bh_{q_i,q_{i+1},\ldots,q_k}(t_i-\tau,t_{i+1},\ldots,t_k) \bu(\tau+T_{i-1}) d\tau,
\end{align}
\normalsize
where the functions $\bg, \bh$ are defined for $k \geqslant 1$, as follows,
\small
\begin{equation}\label{init_state}
\bg_{q_1,q_2,\ldots,q_k}(t_1,t_2,\ldots,t_k) = \bC_{q_k}e^{\tilde{\bA}_{q_k}t_k} \tilde{\bK}_{q_{k-1},q_{k}} e^{\tilde{\bA}_{q_{k-1}}t_{k-1}} \tilde{\bK}_{q_{k-2},q_{k-1}} \cdots \tilde{\bK}_{q_{1},q_{2}} e^{\tilde{\bA}_{q_1}t_1} \bx_0,
\end{equation}
\begin{equation}\label{kern}
\bh_{q_1,q_2,\ldots,q_k}(t_1,t_2,\ldots,t_k) = \bC_{q_k}e^{\tilde{\bA}_{q_k}t_k} \tilde{\bK}_{q_{k-1},q_{k}} e^{\tilde{\bA}_{q_{k-1}}t_{k-1}} \tilde{\bK}_{q_{k-2},q_{k-1}} \cdots \tilde{\bK}_{q_{1},q_{2}} e^{\tilde{\bA}_{q_1}t_1}\tilde{\bB}_{q_1}.
\end{equation}
\normalsize

Note that, for the functions defined in (\ref{init_state}) and (\ref{kern}) we consider the $\bE_{q_i}$ matrices to be incorporated into the $\bA_{q_i}$ and $\bB_{q_i}$ matrices (i.e. $\tilde{\bA}_{q_i} = \bE_{q_i}^{-1} \bA_{q_i}, \ \tilde{\bB}_{q_i} = \bE_{q_i}^{-1} \bB_{q_i}$). Moreover, the transformed coupling matrices are written accordingly $\tilde{\bK}_{q_i,q_{i+1}} = \bE_{q_{i+1}}^{-1} \bK_{q_i,q_{i+1}}$.


By applying the multivariate Laplace transform of the regular kernels in (\ref{kern}), we construct {\it level k generalized transfer functions} of the system $\Si$, as
\begin{equation}\label{gen_trf_lss}
\bH_{q_1,q_2,...,q_k}(s_1,s_2,...,s_k) = \bC_{q_1} \bPhi_{q_1}(s_1) \bK_{q_2,q_1} \bPhi_{q_2}(s_2) \cdots \bK_{q_{k},q_{k-1}}\bPhi_{q_k}(s_k) \bB_{q_k},
\end{equation}
where $\bPhi_q(s) = (s \bE_q - \bA_q)^{-1}$, $q_j \in \{1,2,...,D\}, \ 1 \leqslant j \leqslant k$ and $k \geqslant 3$. Their definition is similar to the ones corresponding to bilinear systems (see \cite{agi16}).

By using their samples, directly construct (reduced) switched models that interpolate the original model, by means of the Loewner framework, as in \cite{gpa17}.

%

\section{Energy Gramians for LSS with two modes}

\subsection{Setup and notations}

For simplicity of the exposition, we first consider the simplified case $D = 2$ (the LSS system switches between two modes only). This situation is encountered in most of the numerical examples in the literature we came across.
Nevertheless, all the theoretical concepts presented in this section can be generalized for a general number of modes denoted with D (as in Section 4, where the main results are directly presented for the general case).
Depending on the values of the switching signal $\sigma(t)$, the original system $\Si$ switches between the following subsystems,
$$
\Si_1: \ \begin{cases}
\bE_{1} \dot{\bx}_1(t) = \bA_{1}\bx_1(t) + \bB_{1}\bu(t), \\
\by(t) = \bC_{1} \bx_1(t) ,
\end{cases} \ \text{or} \ \ \ \Si_2: \ \begin{cases}
\bE_{2} \dot{\bx}_2(t) = \bA_{2}x_2(t) + \bB_{2}\bu(t), \\
\by(t) = \bC_{2} \bx_2(t) ,
\end{cases}
$$
where $\text{dim}(\Si_1) = n_1$ (i.e. $\bx_1 \in \mathbb{R}^{n_1}$ and $\bE_1,\bA_1 \in \mathbb{R}^{n_1 \times n_1}, \bB_1,\bC_1^T \in \mathbb{R}^{n_1}$) and also $\text{dim}(\Si_2) = n_2$ (i.e. $\bx_2 \in \mathbb{R}^{n_2}$ and $\bE_2,\bA_2 \in \mathbb{R}^{n_2 \times n_2}, \bB_2,\bC_2^T \in \mathbb{R}^{n_2}$). Notice that we allow both the two subsystems to be written in descriptor format (having possibly singular E matrix).

Denote, for simplicity, with $\bK_1$ the coupling matrix when switching from mode 1 to mode 2 (instead of $\bK_{1,2}$) and, with $\bK_2$, the coupling matrix when switching from mode 2 to mode 1 (instead of $\bK_{2,1}$) with $\bK_1 \in \mathbb{R}^{n_2 \times n_1}$ and $\bK_2 \in \mathbb{R}^{n_1 \times n_2}$. 


In the following, for the first two levels we present the generalized kernels, which were previously defined in (\ref{kern}), i.e.,
$$
\text{Level 1}:  \begin{cases} \bh_1(t_1) = \bC_1 e^{\bA_1 t_1} \bB_1, \\  \bh_2(t_2) = \bC_2 e^{\bA_2 t_1} \bB_2. \end{cases}, \  \text{Level 2}: \begin{cases} \bh_{1,2}(t_1,t_2) = \bC_1  e^{\bA_{1} t_1} \bK_2 e^{\bA_{2}t_2} \bB_2, \\ \bh_{2,1}(t_1,t_2) = \bC_2  e^{\bA_{2} t_1} \bK_1 e^{\bA_{1}t_2} \bB_1. \end{cases}
$$ 

\begin{definition}
	
	Consider the LSS, $\hat{\Si} = (n_1,n_2, \{(\hat{\bE}_i,\hat{\bA}_i, \hat{\bB}_i, \hat{\bC}_i)\},\{ \hat{\bK}_{i,j}  \},i,j \in \{1,2\}, \bfz)$ and $\bar{\Si} = (n_1,n_2, \{(\bar{\bE}_i,\bar{\bA}_i, \bar{\bB}_i, \bar{\bC}_i)\},\{ \bK_{i,j}  \},,i,j \in \{1,2\}, \bfz)$. These systems are said to be equivalent if there exist non-singular matrices $\bZ_j^L ,\bZ_j^R$ so that
	$$
	\bar{\bE}_j = \bZ_j^L \hat{\bE}_j \bZ_j^R, \ \  \bar{\bA}_j = \bZ_j^L \hat{\bA}_j \bZ_j^R, \ \ \bar{\bB}_j = \bZ_j^L \hat{\bB}_j , \ \ \bar{\bC}_j = \hat{\bC}_j \bZ_j^R, \ \ j \in \{1,2\},
	$$ 
	and also $\bar{\bK}_1 = \bZ_2^L \hat{\bK}_1 \bZ_1^R, \ \bar{\bK}_2 = \bZ_1^L \hat{\bK}_2 \bZ_2^R$. In this configuration, one can easily show that the transfer functions defined above are the same for each LSS and for all sampling points $s_k$.
\end{definition}

Consider a LSS system $\Si$ as described in (\ref{LSS_def0}) with two operational modes, i.e $D=2$ and $\Omega = \{1,2\}$. Consider $\text{dim}(\Si_k) = n_k$ for $k = 1,2$ and let $\bK_1 \in \mathbb{R}^{n_2 \times n_1}$ and $\bK_2 \in \mathbb{R}^{n_1 \times n_2}$ be the coupling matrices.

\begin{definition}
	For $\nu \in \{1,2\}$, let $\Omega^{\nu,+}$ and $\Omega^{+,\nu}$ be the ordered sets containing all tuples that can be constructed with symbols from the alphabet $\Omega = \{1,2\}$ and that \textit{start} (and respectively \textit{end}) with the symbol $\nu$. Also, no two consecutive characters are allowed to be the same. Hence, explicitly write the new introduced sets as follows:
	\begin{eqnarray}
	&& \Omega^{1,+} = \{( 1 ), ( 1,2), (1,2,1), \ldots \}, \ \ \Omega^{2,+} = \{(2),(2,1),(2,1,2),\ldots \}, \\
	&& \Omega^{+,1} = \{ (1),(2,1),(1,2,1),\ldots \}, \ \ \Omega^{+,2} = \{(2),(1,2),(2,1,2),\ldots \}.
	\end{eqnarray}
\end{definition}

\begin{definition}
	Let the $i^{th}$ unit vector of length $k$ be denoted with 
	$$
	\bfe_i=[0,\ldots,1,\ldots,0]^T\in\mathbb{R}^{ k}, \ \ \bfe_i(\ell) = 1, \ \ \text{if} \ \ell = i \ \text{and} \ \bfe_i(\ell) = 0, \ \text{else}.
	$$
	In some contexts we may use the alternative notation $\bfe_{i,k}$ to emphasize its dimension $k$. The identity matrix $\bI_k  \in\mathbb{R}^{ k \times k}
	$ can be written as $\bI_k = [\bfe_{1,k} \ \bfe_{2,k} \ \ldots \ \bfe_{k,k}] $. 
	Also, let $\bfz_{k,\ell} \in \mathbb{R}^{k \times \ell}$ be an all zero matrix. When $k = \ell$, we use the notation $\bfz_k =\in \mathbb{R}^{ k \times k}$ or simply $\bfz$ when the dimension is clearly inferred.
\end{definition}

\subsection{Level k switching - an intermediate step}

\begin{definition}
	In the succeeding sections we analyze LSS $\Si$ for which the $\bE$ matrices corresponding to all of the subsystems $\Si_q$ are identity matrices, i.e. $\bE_q = \bI_{n_q}, \ q \in  \Omega$. Hence we propose an alternative definition for the dynamics of the LSS, i.e.
	\begin{equation}\label{LSS_def}
	\Si: \begin{cases}
	\dot{\bx}(t) = \bA_{\sigma(t)} \bx(t) + \bB_{\sigma(t)}\bu(t), \ \ \bx(t) = \bx_0, \\
	\by(t) = \bC_{\sigma(t)} \bx(t).
	\end{cases}
	\end{equation}
	We again use the compact notation $\Si = (n_1,n_2,\ldots,n_D, \{(\bA_q, \bB_q, \bC_q)| q \in \Omega\},$ $\{ \bK_{q_{i},q_{i+1}} | q_i, q_{i+1} \in \Omega \}, \bx_0)$. The other parameters and notations remain as in (\ref{LSS_def0}).
\end{definition}

\subsubsection{Reachability Gramians}

Introduce the following level $k$ energy functional $\bg^r_{q_1,q_2,\ldots,q_k}(t_1,t_2,\ldots,t_k) : \mathbb{R}^k \rightarrow \mathbb{R}^{m_{q_k}}$, corresponding to the switching sequence $(q_1,q_2,\ldots,q_k) \in \Omega^k$, as
\begin{equation}\label{grqdef}
\bg^r_{q_1,q_2,\ldots,q_k}(t_1,t_2,\ldots,t_k) = e^{\bA_{q_1}t_1} \bK_{q_{2},q_{1}} e^{\bA_{q_{2}}t_{2}} \bK_{q_{3},q_{2}} \cdots \bK_{q_{k},q_{k-1}} e^{\bA_{q_k}t_k} \bB_{q_k}.
\end{equation}

By fixing the first element of the tuple $(q_1,q_2,\ldots,q_k)$, i.e., $q_1 \in \{1,2\}$, note that $(q_1,q_2,\ldots,q_k)$ can either be an element of $\Omega^{1,+}$ or of $ \Omega^{2,+}$ (as introduced in Definition 4).

If we choose $q_1 = 1$, then it follows that $(q_1,q_2,\ldots,q_k) \in \Omega^{1,+}$. Examples of energy functionals associated to sequences from $\Omega^{1,+}$, are for instance the following
\begin{align*}
\bg^r_1(t_1) = e^{\bA_1 t_1} \bB_1, \ \  \bg^r_{1,2}(t_1,t_2) = e^{\bA_1 t_1} \bK_2 e^{\bA_2 t_2}  \bB_2, \\  \bg^r_{1,2,1}(t_1,t_2,t_3) = e^{\bA_1 t_1} \bK_2 e^{\bA_2 t_2} \bK_1 e^{\bA_1 t_3}  \bB_1, \ \ldots
\end{align*}
\noindent
In general, compute the following level $k$ infinite Gramian corresponding to mode $q_1 \in \{1,2\}$  by calculating the inner product of the energy functional associated to the length $k$ switching sequence $(q_1,q_2,\ldots,q_k) \in \Omega^{q_1,+}$ with itself, as
\small
\begin{equation}\label{Pkdef}
\cP_{q_1}^{(k)} = \int_{0}^\infty \cdots \int_{0}^\infty   \bg^r_{q_1,q_2,\ldots,q_k}(t_1,t_2,\ldots,t_k) \big{(}\bg^r_{q_1,q_2,\ldots,q_k}(t_1,t_2,\ldots,t_k) \big{)}^T dt_1 dt_2 \ldots dt_k.
\end{equation}
\normalsize
\noindent
By making use of the recurrence relation
$$
\bg^r_{q_1,q_2,\ldots,q_k}(t_1,t_2,\ldots,t_k) = \big{(} e^{\bA_{q_1}t_1} \bK_{q_{2},q_{1}} \big{)} \bg^r_{q_2,q_3,\ldots,q_k}(t_2,t_3,\ldots,t_k),
$$
it follows that the $k^{\rm th}$ Gramian corresponding to mode 1 (or respectively mode 2) can be written in terms of the $(k-1)^{\rm th}$ Gramian corresponding to mode 2 (or mode 1), as
\begin{align}\label{Pkprop}
\cP_{q_1}^{(k)} &=  \int_{0}^\infty \cdots \int_{0}^\infty   \big{(} e^{\bA_{q_1}t_1} \bK_{q_{2},q_{1}} \big{)} \bg^r_{q_2,\ldots,q_k}(t_2,\ldots,t_k) \big{(} \bg^r_{q_2,\ldots,q_k}(t_2,\ldots,t_k) \big{)}^T  \nonumber \\ & \big{(} e^{\bA_{q_1}t_1} \bK_{q_{2},q_{1}} \big{)}^T dt_1 \ldots dt_k 
= \int_{0}^\infty  e^{\bA_{q_1}t_1} \bK_{q_{2},q_{1}} \Big{(} \int_{0}^\infty  \bg_{q_2,\ldots,q_k}(t_2,\ldots,t_k) \nonumber \\ & \big{(}  \bg^r_{q_2,\ldots,q_k}(t_2,\ldots,t_k) \big{)}^T   dt_2 \ldots dt_k \Big{)}  \bK_{q_{2},q_{1}}^T e^{\bA_{q_1}^Tt_1} dt_1 \nonumber \\
&= \int_{0}^\infty  e^{\bA_{q_1}t_1} \bK_{q_{2},q_{1}}  \cP_{q_2}^{(k-1)} \bK_{q_{2},q_{1}}^T e^{\bA_{q_1}^Tt_1} dt_1.
\end{align}

\begin{proposition}
	Next introduce the linear reachability Gramians for the case with no switching. They are denoted with $\cP_{q}^{(1)}$, corresponding to mode $q \in \{1,2\}$, and can be defined as
	\begin{equation}\label{def_lin_gram_reach}
	\cP_{q}^{(1)} = \int_0^\infty \bg^r_q(t) \big{(}  \bg^r_q(t) \big{)}^T dt = \int_0^\infty e^{\bA_q t} \bB_q \bB_q^T  e^{\bA_q^T t} dt.
	\end{equation}
	It is a well known result that $\cP_{q}^{(1)}$ satisfies the following Lyapunov equation:
	\begin{align}\label{lin_gram_reach}
	\bA_q \cP_q^{(1)} + \cP_q^{(1)} \bA_q^T + \bB_q \bB_q^T = \bfz.
	\end{align}
\end{proposition}

\begin{proposition}
	The level k reachability Gramians corresponding to modes 1 and 2 can be computed by iteratively solving the coupled systems of linear equations:
	\begin{align} 
	\bA_1 \cP_1^{(k)} + \cP_1^{(k)} \bA_1^T + \bK_2 \cP_2^{(k-1)}  \bK_2^T  = \bfz, \label{klevel_reach1} \\ \bA_2 \cP_2^{(k)} + \cP_2^{(k)} \bA_2^T  + \bK_1 \cP_1^{(k-1)}  \bK_1^T = \bfz,
	\label{klevel_reach2}
	\end{align}
	where $k>1$ and the starting point is represented by the linear Gramians (with no switching) $\bP_{q_1}^{(1)}$ in (\ref{lin_gram_reach}) that correspond to the first level.
\end{proposition}
\noindent
{\bf{Proof of Proposition 2}}.
By multiplying the equality in (\ref{Pkprop}) with $\bA_{q_1}$ to the left and with $\bA_{q_1}^T$ to the right, we write
$$
\bA_{q_1} \cP_{q_1}^{(k)} + \cP_{q_1}^{(k)} \bA_{q_1}^T = \int_{0}^\infty \bA_{q_1}  e^{\bA_{q_1}t_1} \bK_{q_{2},q_{1}}  \cP_{q_2}^{(k-1)} \bK_{q_{2},q_{1}}^T e^{\bA_{q_1}^Tt_1} dt_1
$$
$$
+ \int_{0}^\infty  e^{\bA_{q_1}t_1} \bK_{q_{2},q_{1}}  \cP_{q_2}^{(k-1)} \bK_{q_{2},q_{1}}^T e^{\bA_{q_1}^Tt_1} \bA_{q_1}^T  dt_1
$$
$$
= \int_{0}^\infty \frac{d}{dt_1} \Big{(}  e^{\bA_{q_1}t_1} \bK_{q_{2},q_{1}}  \cP_{q_2}^{(k-1)} \bK_{q_{2},q_{1}}^T e^{\bA_{q_1}^Tt_1} dt_1 \Big{)} = - \bK_{q_{2},q_{1}}  \cP_{q_2}^{(k-1)} \bK_{q_{2},q_{1}}^T.
$$
Hence it follows that, for $q_1,q_2 \in \{1,2\}$ with $q_1 \neq q_2$, we write
$$
\bA_{q_1} \cP_{q_1}^{(k)} + \cP_{q_1}^{(k)} \bA_{q_1}^T + \bK_{q_{2},q_{1}}  \cP_{q_2}^{(k-1)} \bK_{q_{2},q_{1}}^T = \bfz,
$$
which proves the statements in (\ref{klevel_reach1}) and (\ref{klevel_reach2}). \sq

\subsubsection{Observability Gramians}

Introduce the following level $k$ energy functional $\bg^o_{q_k,\ldots,q_2,q_1}(t_k,\ldots,t_2,t_1) : \mathbb{R}^k \rightarrow \mathbb{R}^{p_{q_k}}$, corresponding to the switching sequence $(q_k,\ldots,q_2,q_1) \in \Omega^k$, as
\small
\begin{equation}\label{goqdef}
\bg^o_{q_k,q_{k-1},\ldots,q_1}(t_k,\ldots,t_2,t_1) =  \bC_{q_k}e^{\bA_{q_k}t_k} \bK_{q_{k-1},q_{k}} e^{\bA_{q_{k-1}}t_{k-1}} \bK_{q_{k-2},q_{k-1}} \cdots \bK_{q_{1},q_{2}} e^{\bA_{q_1}t_1}
\end{equation}
\normalsize

By fixing the last element of the tuple, i.e., $q_1 \in \{1,2\}$, note that $(q_k,\ldots,q_2,q_1)$ can either be an element of $\Omega^{+,1}$ or of $ \Omega^{+,2}$ (as introduced in Definition 4).

If $q_1 = 1$ is chosen, then it follows that $(q_{k},q_{k-1},\ldots,q_1) \in \Omega^{1,+}$. Examples of energy functionals associated to sequences from $\Omega^{+,1}$, are the following
\begin{align*}
\bg^o_1(t_1) = \bC_1 e^{\bA_1 t_1}, \ \  \bg^o_{2,1}(t_2,t_1) = \bC_2 e^{\bA_2 t_2} \bK_1 e^{\bA_1 t_1} , \\  \bg^o_{1,2,1}(t_3,t_2,t_1) = \bC_1 e^{\bA_1 t_3} \bK_2 e^{\bA_2 t_2} \bK_1 e^{\bA_1 t_1}, \ \ldots
\end{align*}
\noindent
Compute the following level $k$ infinite Gramian corresponding to mode $q_1 \in \{1,2\}$  by calculating the inner product of the energy functional associated to the length $k$ switching sequence $(q_1,q_2,\ldots,q_k) \in \Omega^{q_1,+}$ with itself, as
\begin{equation}\label{Qkdef}
\cQ_{q_1}^{(k)} = \int_{0}^\infty \cdots \int_{0}^\infty  \big{(} \bg^o_{q_k,\ldots,q_2,q_1}(t_k,\ldots,t_2,t_1) \big{)} ^T \bg^o_{q_k,\ldots,q_2,q_1}(t_k,\ldots,t_2,t_1) dt_1\ldots dt_k.
\end{equation}
\noindent
By using the following recurrence relation,
$$
\bg^o_{q_k,\ldots,q_2,q_1}(t_k,\ldots,t_2,t_1) =  \bg^o_{q_k,\ldots,q_3,q_2}(t_k,\ldots,t_3,t_2) \big{(}  \bK_{q_{1},q_{2}} e^{\bA_{q_1}t_1} \big{)},
$$
the $k^{\rm th}$ observability Gramian corresponding to mode 1 (or respectively mode 2) can be written in terms of the $(k-1)^{\rm th}$ observability Gramian corresponding to mode 2 (or respectively mode 1), as
\begin{align}\label{Qkprop}
\cQ_{q_1}^{(k)} &= \int_{0}^\infty \cdots \int_{0}^\infty   \big{(}  \bK_{q_{1},q_{2}} e^{\bA_{q_1}t_1} \big{)}^T \big{(} \bg^o_{q_k,\ldots,q_2}(t_k,\ldots,t_2) \big{)}^T  \bg^o_{q_k,\ldots,q_2}(t_k,\ldots,t_2) \nonumber \\ &   \big{(}  \bK_{q_{1},q_{2}} e^{\bA_{q_1}t_1} \big{)} dt_1 \ldots dt_k  
= \int_{0}^\infty  e^{\bA^T_{q_1}t_1} \bK_{q_{1},q_{2}}^T  \Big{(} \int_{0}^\infty \big{(} \bg^o_{q_k,\ldots,q_2}(t_k,\ldots,t_2) \big{)}^T \nonumber \\ & \bg^o_{q_k,\ldots,q_2}(t_k,\ldots,t_2)  dt_2 \ldots dt_k \Big{)}  \bK_{q_{1},q_{2}} e^{\bA_{q_1}t_1} dt_1 \nonumber \\
&= \int_{0}^\infty  e^{\bA^T_{q_1}t_1} \bK_{q_{1},q_{2}}^T  \cQ_{q_2}^{(k-1)} \bK_{q_{1},q_{2}} e^{\bA_{q_1}t_1} dt_1.
\end{align}

\begin{proposition}
	The linear observability Gramian (for the case with no switching) $\cQ_{q}^{(1)}$ which corresponds to mode $q \in \{1,2\}$, can be written as
	\begin{equation}\label{def_lin_gram_obs}
	\cQ_{q}^{(1)} = \int_0^\infty \big{(} \bg^o_q(t) \big{)}^T \bg^o_q(t) dt = \int_0^\infty e^{\bA_q^T t} \bC_q^T \bC_q  e^{\bA_q t} dt.
	\end{equation}
	It is a well known result that $\cQ_{q}^{(1)}$ satisfies the following Lyapunov equation:
	\begin{align}\label{lin_gram_obs}
	\bA_q^T \cQ_q^{(1)} + \cQ_q^{(1)} \bA_q + \bC_q^T \bC_q = \bfz.
	\end{align}
\end{proposition}

\begin{proposition}
	The level k observability Gramians corresponding to modes 1 and 2 can be computed by iteratively solving the coupled systems of linear equations (for $k>1$)
	\begin{align} 
	\bA_1^T \cQ_1^{(k)} + \cQ_1^{(k)} \bA_1 + \bK_1^T \cQ_2^{(k-1)}  \bK_1  = \bfz, \label{klevel_obs1} \\ \bA_2^T \cQ_2^{(k)} + \cQ_2^{(k)} \bA_2  + \bK_2^T \cQ_1^{(k-1)}  \bK_2 = \bfz,
	\label{klevel_obs2}
	\end{align}
	where the starting point is represented by the linear Gramians (with no switching) $\Omega_{q_1}^{(1)}$ in (\ref{lin_gram_obs}) that correspond to the first level.
\end{proposition}
\noindent
{\bf{Proof of Proposition 4}}.
By multiplying the identity in (\ref{Qkprop}) with $\bA_{q_1}^T$ to the left and with $\bA_{q_1}$ to the right, we write
$$
\bA_{q_1}^T \cQ_{q_1}^{(k)} + \cQ_{q_1}^{(k)} \bA_{q_1}= \int_{0}^\infty \bA_{q_1}^T  e^{\bA^T_{q_1}t_1} \bK_{q_{1},q_{2}}^T  \cQ_{q_2}^{(k-1)} \bK_{q_{1},q_{2}} e^{\bA_{q_1}t_1} dt_1
$$
$$
+ \int_{0}^\infty  e^{\bA^T_{q_1}t_1} \bK^T_{q_{2},q_{1}}  \cQ_{q_2}^{(k-1)} \bK_{q_{1},q_{2}} e^{\bA_{q_1}t_1} \bA_{q_1}  dt_1
$$
$$
= \int_{0}^\infty \frac{d}{dt_1} \Big{(}  e^{\bA^T_{q_1}t_1} \bK_{q_{1},q_{2}}^T  \cQ_{q_2}^{(k-1)} \bK_{q_{1},q_{2}} e^{\bA_{q_1}t_1} dt_1 \Big{)} = - \bK_{q_{1},q_{2}}^T  \cQ_{q_2}^{(k-1)} \bK_{q_{2},q_{1}}.
$$
Hence it follows that, for $q_1,q_2 \in \{1,2\}$ with $q_1 \neq q_2$, we write
$$
\bA_{q_1}^T \cQ_{q_1}^{(k)} + \cQ_{q_1}^{(k)} \bA_{q_1} + \bK_{q_{1},q_{2}}^T  \cQ_{q_2}^{(k-1)} \bK_{q_{1},q_{2}} = \bfz.
$$
which proves the statements in (\ref{klevel_obs1}) and (\ref{klevel_obs2}). \sq

\subsection{Infinite Gramians and Lyapunov equations}

\begin{definition}
	Introduce the infinite reachability Gramian $\cP_{q_1}$ corresponding to mode $q_1 \in \{1,2\}$ of the LSS system $\Si$ as
	\small
	\begin{align*}
	\cP_{q_1} = \sum_{k=1}^\infty \int_{0}^\infty \cdots \int_{0}^\infty   \bg^r_{q_1,q_2,\ldots,q_k}(t_1,t_2,\ldots,t_k)  \big{(}\bg^r_{q_1,q_2,\ldots,q_k}(t_1,t_2,\ldots,t_k) \big{)}^T  dt_1  \ldots dt_k, 
	\end{align*}
	\normalsize
	\vspace{-4mm}
	\begin{align}\label{defPgram}
	\Rightarrow \cP_{q_1} =  \sum_{k=1}^\infty \cP_{q_1}^{(k)} = \cP_{q_1}^{(1)} + \cP_{q_1}^{(2)} + \ldots, 
	\end{align}
	in terms of the multivariate functions $\bg_q^r$ in (\ref{grqdef}) or matrices $\cP_{q_1}^{(k)}$ in (\ref{Pkdef}).
\end{definition}

Note that $\cP_{q_1}$ is computed by taking into account the inner products of energy functionals associated to all possible switching sequences (of any length $k$) that start in mode $q_1$.

\begin{proposition}
	The infinite reachability Gramians defined in (\ref{defPgram}), satisfy the following system of generalizaed coupled Lyapunov equations
	\begin{equation}\label{PgramLyap}
	\begin{cases} \bA_1 \cP_1 +\cP_1 \bA_1^T+ \bK_2 \cP_2 \bK_2^T+ \bB_1 \bB_1^T = \bfz, \\ \bA_2 \cP_2 +\cP_2 \bA_2^T+ \bK_1 \cP_1 \bK_1^T+ \bB_2 \bB_2^T = \bfz. \end{cases}
	\end{equation}
\end{proposition}
{\bf{Proof of Proposition 5}}.
By adding the equalities stated in (\ref{klevel_reach1}) and (\ref{klevel_reach2}) for $k>2$ as well as the one corresponding to $k=1$ (in (\ref{def_lin_gram_reach})), it 
follows that
\small
$$
\big{(} \bA_{q_1} \cP_{q_1}^{(1)} + \cP_{q_1}^{(1)} \bA_{q_1}^T + \bB_{q_1}\bB_{q_1}^T \big{)} + \sum_{k=2}^ \infty \big{(} \bA_{q_1} \cP_{q_1}^{(k)} + \cP_{q_1} \bA_{q_1}^T + \bK_{q_2,q_1} \cP_{q_2}^{(k-1)}  \bK_{q_2,q_1}^T \big{)} = \bfz \\
$$
\normalsize
\vspace{-4mm}
\begin{align*}
&  \Rightarrow \bA_{q_1} \big{(} \sum_{k=1}^ \infty \cP_{q_1}^{(k)} \big{)} + \big{(} \sum_{k=1}^ \infty \cP_{q_1}^{(k)} \big{)} \bA_{q_1}^T +  \bK_{q_2,q_1} \big{(} \sum_{k=1}^ \infty  \cP_{q_1}^{(k)} \big{)} \bK_{q_2,q_1}^T+ \bB_{q_1}\bB_{q_1}^T = \bfz.
\\
&  \Rightarrow \bA_{q_1} \cP_{q_1} + \cP_{q_1} \bA_{q_1}^T +  \bK_{q_2,q_1} \cP_{q_1} \bK_{q_2,q_1}^T+ \bB_{q_1}\bB_{q_1}^T = \bfz \ \ \forall q_1 \neq q_2 \in \{1,2\},
\end{align*}
which shows the validity of the equalities introduced in (\ref{PgramLyap}).  \sq

\begin{remark}
	Write the  matrices $\{\cP_q, \bA_q,\bB_q,\bC_q\}$, $q \in \{1,2\}$ and $\{\bK_{q_1,q_2}\}$, $q_1, q_2 \in \{1,2\}$  in block-diagonal format, as
	\begin{equation}\label{notationD}
	\bX_{\bD} =\left[ \begin{array}{cc}
	\bX_1 & \bfz \\ \bfz & \bX_2
	\end{array} \right], \ \bX \in \{\bA,\bB,\bC, \cP\}, \ \ \ \bK_{\scriptsize\reflectbox{\bD}} =\left[ \begin{array}{cc}
	\bfz & \bK_1  \\  \bK_2 & \bfz
	\end{array} \right].
	\end{equation} 
	Hence, instead of solving the two equations in (\ref{QgramLyap}) separately, one can solve one equation
	\begin{equation}\label{PgramDLyap}
	\bA_{\bD} \bP_\bD + \bP_\bD \bA_\bD^T +\bK_{\scriptsize\reflectbox{\bD}} \bP_\bD \bK_{\scriptsize\reflectbox{\bD}}^T+ \bB_\bD \bB_\bD^T = \bfz,
	\end{equation}
	and recover the reachability Gramians $\cP_1$ and $\cP_2$ as block diagonal entries of $\bP_{\bD}$.
\end{remark}

\begin{definition}
	Introduce the infinite observability Gramian $\cQ_{q_1}$ corresponding to mode $q_1 \in \{1,2\}$ of the LSS system $\Si$ as
	\small
	\begin{equation*}
	\cQ_{q_1} = \sum_{k=1}^\infty \int_{0}^\infty \cdots \int_{0}^\infty \big{(}  \bg^o_{q_k,\ldots,q_2,q_1}(t_k,\ldots,t_2,t_1) \big{)}^T \bg^o_{q_k,\ldots,q_2,q_1}(t_k,\ldots,t_2,t_1) \ dt_1  dt_2 \ldots dt_k 
	\end{equation*}
	\normalsize
	\vspace{-4mm}
	\begin{equation} \label{defQgram}
	\cQ_{q_1} = \sum_{k=1}^\infty \cQ_{q_1}^{(k)} = \cQ_{q_1}^{(1)} + \cQ_{q_1}^{(2)} + \ldots
	\end{equation}
\end{definition}
\noindent
Note that $\cQ_{q_1}$ is computed by taking into account the inner products of  energy functionals associated to all possible switching sequences (of any length $k$) that end in mode $q_1$.

\begin{proposition}
	The infinite reachability Gramians defined in (\ref{defQgram}), satisfy the following system of generalizaed coupled Lyapunov equations
	\begin{equation}\label{QgramLyap}
	\begin{cases} \bA_1^T \cQ_1 +\cQ_1 \bA_1+ \bK_1^T \cQ_2 \bK_1+ \bC_1^T \bC_1 = \bfz \\ \bA_2^T \cQ_2 +\cQ_2 \bA_2+ \bK_2^T \cQ_1 \bK_2+ \bC_2^T \bC_2 = \bfz \end{cases}
	\end{equation}
	in terms of the multivariate functions $\bg_q^o$ in (\ref{grqdef}) and matrices $\cQ_{q_1}^{(k)}$ in (\ref{Qkdef}).
\end{proposition}
{\bf{Proof of Proposition 6}}.
By adding the equalities stated in (\ref{klevel_obs1}) and (\ref{klevel_obs2}) for $k>2$ as well as the one corresponding to $k=1$ (in (\ref{def_lin_gram_obs})), it follows that
\begin{align*}
& \big{(} \bA_{q_1}^T \cQ_{q_1}^{(1)} + \cQ_{q_1}^{(1)} \bA_{q_1} + \bC_{q_1}^T \bC_{q_1}^T \big{)} + \sum_{k=2}^ \infty \big{(} \bA_{q_1}^T \cQ_{q_1}^{(k)} + \cQ_{q_1} \bA_{q_1} + \bK_{q_1,q_2}^T \cQ_{q_2}^{(k-1)}  \bK_{q_1,q_2} \big{)} = \bfz \\
&  \Rightarrow \bA_{q_1}^T \big{(} \sum_{k=1}^ \infty \cQ_{q_1}^{(k)} \big{)} + \big{(} \sum_{k=1}^ \infty \cQ_{q_1}^{(k)} \big{)} \bA_{q_1} +  \bK_{q_1,q_2}^T \big{(} \sum_{k=1}^ \infty  \cQ_{q_1}^{(k)} \big{)} \bK_{q_1,q_2}+ \bC_{q_1}^T \bC_{q_1} = \bfz.
\\
&  \Rightarrow \bA_{q_1}^T \cQ_{q_1} + \cQ_{q_1} \bA_{q_1} +  \bK_{q_1,q_2}^T \cQ_{q_1} \bK_{q_1,q_2}+ \bC_{q_1}^T \bC_{q_1} = \bfz \ \ \forall q_1 \neq q_2 \in \{1,2\},
\end{align*}
which shows the validity of the equalities presented in (\ref{QgramLyap}). \sq

\begin{remark}
	{\rm
		Additional to (\ref{notationD}), write the  matrices $\{\cQ_q\}, q \in \{1,2\}$  in block-diagonal format, as $\bQ_\bD =\left[ \begin{array}{cc}
		\cQ_1 & \bfz \\ \bfz & \cQ_2
		\end{array} \right]$. Hence, instead of solving the two equations in (\ref{QgramLyap}) separately, one can solve one equation
		\begin{equation}\label{QgramDLyap}
		\bA_{\bD}^T \bQ_\bD + \bQ_\bD \bA_\bD +\bK_{\scriptsize\reflectbox{\bD}}^T \bQ_\bD \bK_{\scriptsize\reflectbox{\bD}}+ \bC_\bD^T \bC_\bD = \bfz,
		\end{equation}
		and recover the observability Gramians  as the block diagonal entries of $\bQ_{\bD}$.
	}
\end{remark}

\begin{definition}\label{def1}
	We assume both $\bA_1$ and $\bA_2$ matrices have eigenvalues with negative real part, i.e. $\text{Re}(\lambda_i (\bA_k)) <0, \ k \in \{1,2\}$. Hence, the same property applies for $\bA_\bD$. The system $\dot{\bx} = \bA_\bD \bx$ is asymptotically stable , or in short, $\bA_\bD$ is stable, if there exist real scalars $\beta > 0$ and $0 < \alpha  \leqslant -\max_i (\text{Re}(\lambda_i (\bA_\bD)))$, such that:
	$$
	\Vert e^{\bA_\bD t} \Vert \leqslant \beta e^{-\alpha t}.
	$$
\end{definition}

The following result from \cite{zl02} addresses the existence of the new defined Gramians. In a nutshell, it states that this holds if the norm of the coupling matrices is sufficiently small. 

\begin{proposition}
	The reachability and observability Gramians  in (\ref{defPgram}), (\ref{defQgram}) exist if
	\begin{equation}\label{existPQ}
	\bA_\bD \  \text{is  stable and} \ \ \Vert \bK_\bD \Vert = \max(\Vert \bK_1 \Vert, \Vert \bK_2 \Vert ) \leqslant \frac{\sqrt{2 \alpha}}{\beta}.
	\end{equation}
\end{proposition}

For high order examples, it is not trivial to solve such generalized Lyapunov equations as (\ref{PgramDLyap}) and (\ref{QgramDLyap}). A possible approach is to approximate these solutions with truncated sums of positive definite matrices,
\begin{equation}
\bP_\bD \approx \sum_{k=1}^H \bP_\bD^{(k)}, \ \ \bQ_\bD \approx \sum_{k=1}^H \bQ_\bD^{(k)}, \ \ H \geqslant 1,
\end{equation}
where $\bP_\bD^{(k)}$ and $\bQ_\bD^{(k)}$ can be written as solutions of regular Lyapunov equations, 
\begin{align*}
& \bA_{\bD} \bP_\bD^{(k)} + \bP_\bD^{(k)} \bA_\bD +\bK_{\scriptsize\reflectbox{\bD}} \bP_\bD^{(k-1)} \bK_{\scriptsize\reflectbox{\bD}} ^T = \bfz, \\
& \bA_{\bD}^T \bQ_\bD^{(k)} + \bQ_\bD^{(k)} \bA_\bD +\bK_{\scriptsize\reflectbox{\bD}}^T \bQ_\bD^{(k-1)} \bK_{\scriptsize\reflectbox{\bD}} = \bfz, \ \ k \geqslant 2.
\end{align*}
For practical applications, solving many such Lyapunov equations is expensive. One can compute low rank factors instead of the full solutions to speed up the calculations ad avoid memory problems (for example, by using the toolbox in \cite{skb16}).

\subsection{Extension to LSS with D modes}

Let $\Omega = \{1,2,\ldots,D\}, \ D \geqslant 2$ and fix the starting mode $q_1 \in \Omega$. Introduce the switching scenario $(q_1,q_2,\ldots,q_k) \in \Omega^k$. Since we exclude equal neighboring modes, i.e.  $q_j \neq q_{j+1}, \ 1 \leqslant j \leqslant k-1$, it follows that there are $(D-1)^{k-1}$ ways of choosing such a switching sequence $(q_1,q_2,\ldots,q_k)$. For $D=2$, there was only one possible sequence chosen uniquely. 

For general number of modes $D$, we have to take into consideration the inner products corresponding to all sequences; hence adapt the definition of $\cP_{q_1}^{(k)}$ from (\ref{Pkdef}) as follows
\begin{align}\label{Pkdef_gen}
\cP_{q_1}^{(k)} &= \int_{0}^\infty \cdots \int_{0}^\infty  \sum_{q_2 = 1, \ q_2 \neq q_1}^D \ldots \sum_{q_k = 1, \ q_k \neq q_{k-1}}^D   \bg^r_{q_1,q_2,\ldots,q_k}(t_1,t_2,\ldots,t_k) \nonumber \\ & \big{(}\bg^r_{q_1,q_2,\ldots,q_k}(t_1,t_2,\ldots,t_k) \big{)}^T dt_1 dt_2 \ldots dt_k.
\end{align}
Again, one can write a recurrence relation by fixing the mode indexes $q_3,\ldots, q_k$,
$$
\bg^r_{q_1,q_2,\ldots,q_k}(t_1,t_2,\ldots,t_k) = \sum_{q_2 = 1, \ q_2 \neq q_1}^D \big{(} e^{\bA_{q_1}t_1} \bK_{q_{2},q_{1}} \big{)} \bg^r_{q_2,q_3,\ldots,q_k}(t_2,t_3,\ldots,t_k).
$$
Next, it follows that the $k^{\rm th}$ reachability Gramian corresponding to mode $q_1$ can be written in terms of the $(k-1)^{\rm th}$ reachability Gramians corresponding to modes $\Omega \setminus \{q_1\}$, as
\begin{align}\label{Pkprop_gen}
\cP_{q_1}^{(k)} =  \int_{0}^\infty  \sum_{q_2 = 1, \ q_2 \neq q_1}^D e^{\bA_{q_1}t_1} \bK_{q_{2},q_{1}}  \cP_{q_2}^{(k-1)} \bK_{q_{2},q_{1}}^T e^{\bA_{q_1}^Tt_1} dt_1.
\end{align}


\begin{definition}
	\rm When considering the general case with $D \geqslant 2$ switching modes, define the infinite reachability Gramian corresponding to mode $q_1 \in \Omega$, as 
	\begin{equation}\label{defP_gen}
	\cP_{q_1} = \sum_{q_1 = 1}^D \cP_{q_1}^{(k)}.
	\end{equation}
	Moreover, the equations satisfied by the reachability Gramians $\cP_i$, for $i \in \{1,2,\ldots,D\}$ can be extended from (\ref{PgramLyap}), as follows
	\begin{equation}\label{PgramLyap_gen}
	\bA_i \cP_i +\cP_i \bA_i^T+ \sum_{j=1, \ j \neq i}^D \bK_{j,i} \cP_j \bK_{j,i}^T+ \bB_i \bB_i^T = \bfz,
	\end{equation}
\end{definition}

\begin{definition}
	Similarly, we can write the observability Gramians as,
	\begin{equation}\label{defQ_gen}
	\cQ_{q_1} = \sum_{q_1 = 1}^D \cQ_{q_1}^{(k)}.
	\end{equation}
	Again the  the system of generalized Lyapunov equations 
	\begin{equation}\label{QgramLyap_gen}
	\bA_i^T \cQ_i +\cQ_i \bA_i+ \sum_{j=1, \ j \neq i}^D \bK_{i,j}^T \cQ_j \bK_{i,j}+ \bC_i^T \bC_i = \bfz.
	\end{equation}
	is satisfied by the matrices $\cQ_{i}, \ \in \Omega$.
\end{definition}

The Gramians introduced in Definition 10 and 11 are mainly going to be used for the original possibly large-scale system. In this case, we would like to  avoid computing the Gramians as solutions of LMI's (as in \cite{pwl13}). Additionally, we present a more relaxed definition of Gramians which will turn out to be useful for the reduced low order case. 

\begin{definition}
	\rm The relaxed reachability Gramians $\cP_i > \bfz$   are defined as solutions of the following collection of LMI, for $i \in \{1,2,\ldots,D\}$ and scalar $M>0$
	\begin{equation}\label{PgramLyap_gen}
	\bA_i \cP_i +\cP_i \bA_i^T+ M \cP_i + \bB_i \bB_i^T < \bfz.
	\end{equation}
	Similarly, the relaxed observability Gramians  $\cQ_i$, satisfy  the inequalities for $i \in \{1,2,\ldots,D\}$, 
	\begin{equation}\label{QgramLyap_gen}
	\bA_i^T \cQ_i +\cQ_i \bA_i+ M \cQ_i+ \bC_i^T \bC_i < \bfz.
	\end{equation}
\end{definition}
\noindent

Note that a Gramian is also a relaxed Gramian but the converse is not necessarily true. Next, we will generalize the results form Remark 3 and 4 for the case with D modes.

Let $\tau_k^n:\{1,\ldots,n\} \rightarrow \{1,\ldots,n\}$ be a cyclic permutation of index k where $k \in \{0,1,\ldots,n-1\}$. The explicit rule is given by $\tau_k^n(\ell) = {\rm \overline{mod}}(k+\ell,n), \ \ \ell \in \{1,\ldots,n\}$, while the permutation $\tau_k^n$ can also be written as,
\begin{equation}
\tau_k^n = \left( \begin{array}{cccc} 1 & 2 & \ldots & n \\  {\rm \overline{mod}}(k+1,n) \ \ & {\rm \overline{mod}}(k+2,n) \ \ & \ldots & \ \ {\rm \overline{mod}}(k+n,n)
\end{array} \right),
\end{equation}
${\rm \overline{mod}}: \{1,\ldots,2n-1\} \rightarrow \{1,\ldots,n\}, \ \ {\rm \overline{mod}}(k,n) = \begin{cases} k, \ \text{if} \ 1 \leqslant k \leqslant n-1 \\   n, \ \text{if} \ k = n \\ k-n, \ \text{if} \ n+1 \leqslant k \leqslant 2n-1 \end{cases}$. 
\noindent
Introduce the permutation matrix $\bPhi_k^n \in \mathbb{R}^{n \times n}$ corresponding to $\tau_k^n$, that has the $\ell^{\rm th}$ row equal to the unit vector $\bfe_{\tau_k^n(\ell),n}^T$. Note that $\bPhi_k^n \bPhi_{n-k}^n = \bI_n$ and $(\bPhi_k^n)^T = \bPhi_{n-k}^n , \ k \in \{0,\ldots,n\}$. For example, write 
$$
\tau_0^3 = \left( \begin{array}{ccc} 1 & 2 & 3 \\ 1 & 2 & 3 \end{array} \right), \ \ \tau_1^3 = \left( \begin{array}{ccc} 1 & 2 & 3 \\ 2 & 3 & 1\end{array} \right), \ \ \tau_2^3 = \left( \begin{array}{ccc} 1 & 2 & 3 \\ 3 & 1 & 2\end{array} \right), \ \text{and} \ \bPhi_2^3 = \left[\begin{array}{ccc} 0 & 0 & 1\\ 1 & 0 & 0\\ 0 & 1 & 0 \end{array}\right].
$$
\begin{remark}
	One can rewrite the $D$ equations stated in (\ref{PgramLyap_gen}) as one equation in the following way,
	\begin{equation}\label{PgramLyap_genD}
	\bA_{\bD} \bP_\bD + \bP_\bD \bA_\bD^T +\sum_{k=1}^{D-1} \bK_{\scriptsize\reflectbox{\bD}_k} \bP_\bD \bK_{\scriptsize\reflectbox{\bD}_k} ^T + \bB_\bD \bB_\bD^T = \bfz.
	\end{equation}
\end{remark}
For all $\bX \in \{\bA,\bB,\bC,\cP,\cQ\}$ and $k \in \{1,\ldots,D-1\}$, consider the notations
\begin{equation}\label{notationDgen}
\bX_{\bD} = \left[ \begin{array}{cccc} \bX_1 & 0 & \ldots & 0\\ 0 & \bX_1 & \ldots & 0\\ 0 & 0 & \ddots & 0  \\ 0 & 0 & \ldots & \bX_D\end{array}\right], \  \bK_{\scriptsize\reflectbox{\bD}_k} = \tilde{\bPhi}_k^{D} \left[ \begin{array}{cccc} \bK_{1,\tau_{D-k}^D(1)} & 0 & \ldots & 0\\ 0 & \bK_{2,\tau_{D-k}^D(2)} & \ldots & 0\\ 0 & 0 & \ddots & 0  \\ 0 & 0 & \ldots & \bK_{D,\tau_{D-k}^D(D)}\end{array}\right],
\end{equation}
where $\tilde{\bPhi}_k^{D} \in \mathbb{R}^{\sum_{i=1}^D n_i \times \sum_{i=1}^D n_i}$ is a block-permutation matrix written in terms of $\bPhi_k^{D}$, by replacing its one entries with  identity matrices $\bI_{n_i}$ of appropriate dimensions. For example, choose $D=3$ and $k=2$, and write the matrix $\tilde{\bPhi}_2^3$ as:
$$
\tilde{\bPhi}_2^3 = \left[\begin{array}{ccc} 0 & 0 & \bI_{n_2}\\ \bI_{n_3} & 0 & 0\\ 0 & \bI_{n_1} & 0 \end{array}\right] \in \mathbb{R}^{(n_1+n_2+n_3) \times (n_1+n_2+n_3)}.
$$
Note that, following the definition of $\bPhi_k^n$, we can write that $\tilde{\bPhi}_2^3 = \bPhi_{n_1+n_3}^{n_1+n_2+n_3}$.
\begin{remark}
	Similarly, we can rewrite the equations in (\ref{QgramLyap_gen}) as only one equation,
	\begin{equation}\label{QgramLyap_genD}
	\bA_{\bD}^T \bQ_\bD + \bQ_\bD \bA_\bD +\sum_{k=1}^{D-1} \bK_{\scriptsize\reflectbox{\bD}_k}^T \bQ_\bD \bK_{\scriptsize\reflectbox{\bD}_k}+ \bC_\bD^T \bC_\bD = \bfz.
	\end{equation}
\end{remark}

\section{Main results}

In this section, we will provide a collection of results that involve the new defined infinite Gramians. In particular, these results will correspond to the more general case with D discrete modes, as presented in Definition 10,11 and 12.

\subsection{Energy bounds relating the input or output signals}

First, we present  the system theoretic interpretation approach; one can write upper and lower bounds of the energy of observation and respectively, of the energy of control in terms of the quantities $\cQ_i$ and $\cP_i$.

\subsubsection{Observability Gramians}

\begin{assumption}
	{ \rm 
		By considering that $\sum_{j=1, \ j \neq i}^{D} \bK_{i,j}^T \cQ_j \bK_{i,j} >0, \ \forall i \in \Omega$, one can show that there exist scalars $M_i >0$ so that
		\begin{equation} \label{condition1Q}
		\sum_{j=1, \ j \neq i}^{D} \bK_{i,j}^T \cQ_j \bK_{i,j} \geqslant M_i \cQ_i, \ \forall i,j \in \Omega.
		\end{equation}
		Additionally, one can also find scalars $\gamma_{i,j} >0$ to satisfy the following inequalities
		\begin{equation} \label{condition2Q}
		\gamma_{i,j} \bK^T_{i,j} \cQ_j \bK_{i,j} < \cQ_i.
		\end{equation}
	}
\end{assumption}

\begin{lemma}
	Given an LSS $\Si$ as defined in (\ref{LSS_def}), consider that the equations in (\ref{QgramLyap_gen}) have positive definite solutions $\cQ_q > \bfz, \ q \in \Omega$. Then, there exists a dwell time $\mu>0$ so that for any switching signal in  (\ref{switch_signal}), with $t_i \geqslant \mu, \ \forall i \geqslant 1$, and zero input $\bu(t) = \bfz$, the following holds
	\begin{equation}\label{energy_bound_Q}
	\bx(0)^T \bQ_{q_1} \bx(0) \geqslant \int_0^t \by^T(s) \by(s) ds, \ \ \forall t>0,
	\end{equation}
	where $q_1 \in \Omega$ represents the index of the first discrete mode in which $\Si$ operates.
\end{lemma}
\noindent
{\bf{Proof of Lemma 1}}. \ 
Consider that the conditions stated in Assumption 41 hold. Introduce $\gamma = \min\limits_{i,j \in \Omega, \ i \neq j}\gamma_{i,j}$ and $M =  \min\limits_{i \in \Omega}M_i$. Choose the minimal dwell times as $\mu = - \frac{\ln \gamma}{M}$. For any piecewise continuous switching signal $\sigma: \mathbb{R} \rightarrow \Omega$ satisfying the conditions in (\ref{switch_signal}) and with minimal dwell time $\mu$, we will prove the bound stated in (\ref{energy_bound_Q}). From (\ref{QgramLyap_gen}) and (\ref{condition1Q}), it follows that
$\bA_i^T \cQ_i +\cQ_i \bA_i+ M_i \cQ_i + \bC_i^T \bC_i \leqslant \bfz$ and furthermore,
\begin{equation}\label{ineqQ}
\bA_i^T \cQ_i +\cQ_i \bA_i+ M \cQ_i + \bC_i^T \bC_i \leqslant \bfz.
\end{equation}
Let $\bx(t)$ the corresponding solution to (\ref{LSS_def}), and  also introduce the functions  $V,W:\mathbb{R}^{n_i} \rightarrow \mathbb{R}$ as
\begin{equation}\label{Vdef}
V(\bx(t)) = \begin{cases} \bx^T(t) \cQ_{q_1} \bx(t), \  t \in [0,t_1] \\ \bx^T(t) \cQ_{q_i} \bx(t),\  t \in (T_{i-1}, T_i], \ i \geqslant 2 \end{cases},
\end{equation}
\begin{equation}\label{Wdef}
W(\bx(t)) = \begin{cases} e^{Mt} \bx(t)^T \cQ_{q_1} \bx(t), \  t \in [0,t_1] \\ e^{M(t-T_{i-1})}V(\bx(t)),  \  t \in (T_{i-1},T_i], \ i \geqslant 2 \end{cases},
\end{equation}
where $T_i = \sum_{\ell=1}^{i} t_\ell$. By considering the uncontrolled case, the input function is considered to be $\bu(t) = 0, \ \forall t$. Using that $\frac{d \bx(t)}{dt} = \bA_{q_i} \bx(t)$, write the derivative of $V(t)$ from (\ref{Vdef}) for $ t \in (T_{i-1},T_i]$,
\begin{equation*}
\frac{\partial V(\bx(t))}{\partial t} = \frac{d \bx^T(t)}{d t} \cQ_{q_i} \bx(t) +  \bx^T(t) \cQ_{q_i}\frac{d \bx(t)}{d t} = \bx^T(t) \big{(} \bA_{q_i}^T \cQ_{q_i} + \cQ_{q_i} \bA_{q_i}  \big{)} \bx(t).
\end{equation*}
For $t \in (T_{i-1},T_i]$, compute the time derivative of $W(\bx(t))$ as defined in (\ref{LSS_def2}) in terms of the one corresponding to $V(x(t))$, as
\begin{align}
\frac{\partial W(\bx(t))}{\partial t} &= Me^{M(t-T_{i-1})}V(\bx(t))+e^{M(t-T_{i-1})} \frac{\partial V(\bx(t))}{\partial t} \nonumber \\ &=
e^{M(t-T_{i-1})} \Big{(} M V(\bx(t))+ \bx^T(t) \big{(} \bA_{q_i}^T \cQ_{q_i} + \cQ_{q_i} \bA_{q_i}  \big{)} \bx(t)  \Big{)} \nonumber \\ &= e^{M(t-T_{i-1})}  \bx^T(t) \big{(} \bA_{q_i}^T \cQ_{q_i} + \cQ_{q_i} \bA_{q_i} + M \cQ_i   \big{)}  \bx(t).  \label{derWx}
\end{align}
By substituting the inequality in (\ref{ineqQ}) into the above relation (\ref{derWx}), and using that $\by(t) = \bC_i \bx(t), \ t \in (T_{i-1},T_i]$, it follows that
\begin{equation}\label{Wderiv}
\frac{\partial W(\bx(t))}{\partial t} \leqslant - e^{M(t-T_{i-1})} \by(t)^T \by(t).
\end{equation}
Introduce the following notation 
\begin{equation}\label{notationQ}
\bx(T_i^+) = \lim\limits_{t \searrow T_i} \bx(t),   \ \ V(\bx(T_i^+)) = \lim\limits_{t \searrow T_i} V(\bx(t)), \  \ W(\bx(T_i^+)) = \lim\limits_{t \searrow T_i} W(\bx(t)).
\end{equation}
By integrating the inequality (\ref{Wderiv}) from $T_{i-1}$ to $t \in (T_{i-1},T_i] $, it follows that
\begin{equation}\label{Wdiff}
W(\bx(t))-W(\bx(T_{i-1}^+)) \leqslant -  \int\limits_{T_{i-1}}^t e^{M(s-T_{i-1})} \by(s)^T \by(s) ds \leqslant -  \int\limits_{T_{i-1}}^t  \by(s)^T \by(s) ds.
\end{equation}
From (\ref{Vdef}) and (\ref{Wdef}), it follows that 
\begin{align}\label{WTi}
W(\bx(T_i)) = e^{M(T_i-T_{i-1})} V(\bx(T_i)) = e^{M t_i} V(\bx(T_i)),
\end{align}
and additionally, using that $\bx(T_i^+) = \bK_{q_i,q_{i+1}} \bx(T_i)$, write
\begin{equation}\label{WTiplus}
W(\bx(T_i^+)) =  V(\bx(T_i^+)) =  \bx^T(T_i) \bK_{q_i,q_{i+1}}^T \cQ_{q_{i+1}} \bK_{q_i,q_{i+1}} \bx(T_i).
\end{equation}
From $(\ref{condition2Q})$ and $(\ref{WTiplus})$ and using that 
$\gamma = \min\limits_{i,j \in \Omega, \ i \neq j}\gamma_{i,j}$, write
\begin{equation}\label{WTiplus2}
W(\bx(T_i^+)) \leqslant \frac{1}{\gamma} \bx(T_i)^T \cQ_i \bx(T_i) = \frac{1}{\gamma} V(\bx(T_i)).
\end{equation}
By combining (\ref{WTi}) and (\ref{WTiplus2}), we can write
\begin{equation}\label{WTiplus3}
W(\bx(T_i^+)) \leqslant \frac{e^{-Mt_i}}{\gamma} W(\bx(T_i)).
\end{equation}
Since switching signals $\sigma$ with minimal dwell time $\mu$ are considered, it follows that $t_i \geqslant \mu \Rightarrow \frac{e^{-Mt_i}}{\gamma} \leqslant \frac{e^{-M \mu}}{\gamma}$. Since, by definition $\mu = -\frac{\ln{\gamma}}{M}$, we get that $\frac{e^{-Mt_i}}{\gamma} \leqslant 1$. Therefore, from (\ref{WTiplus3}), write
\begin{equation}\label{WTiplus4}
W(\bx(T_i^+)) \leqslant W(\bx(T_i)).
\end{equation}
Putting together the inequalities in (\ref{Wdiff}) and (\ref{WTiplus4}), it follows that
\begin{equation}\label{WTiplus5}
W(\bx(T_i)) - W(\bx(T_{i-1})) \leqslant  -  \int\limits_{T_{i-1}}^{T_i}  \by(s)^T \by(s) ds.
\end{equation}
Now using the convention $T_0 = 0$ and adding all the inequalities in (\ref{WTiplus5}), we obtain
\begin{align}\label{WTiplus6}
\sum\limits_{i=1}^\ell W(\bx(T_i)) - W(\bx(T_{i-1})) \leqslant  -  \sum\limits_{i=1}^\ell \int\limits_{T_{i-1}}^{T_i}  \by(s)^T \by(s) ds \nonumber \\ \Rightarrow W(\bx(T_\ell)) - W(\bx(0)) \leqslant - \int\limits_{0}^{T_\ell}  \by(s)^T \by(s) ds.
\end{align}
Since $W(\bx(T_\ell)) = e^{M t_\ell} \bx^T(T_\ell) \cQ_{q_\ell} \bx(T_\ell)  \geqslant 0$, from (\ref{WTiplus6}) it follows that, 
\begin{equation}\label{WTiplus7}
W(\bx(0))   \geqslant  \int\limits_{0}^{T_\ell}  \by(s)^T \by(s) ds, \ \ \forall \ \ell \leqslant 0.
\end{equation}
Now using that $W(\bx(0)) = \bx(0)^T \cQ_{q_1} \bx(0)$, the result in (\ref{energy_bound_Q}) is hence proven.

\subsubsection{Reachability Gramians}

\begin{assumption}
	Considering that $\sum_{j=1, \ j \neq i}^{D} \bK_{j,i} \cP_j \bK_{j,i}^T >0, \ \forall i \in \Omega$,  one can always find scalars $M_i > 0$ such that
	\begin{equation} \label{condition1P}
	\sum_{j=1, \ j \neq i}^{D} \bK_{j,i} \cP_j \bK_{j,i}^T \geqslant M_i \cP_i, \ \forall i \in \Omega.
	\end{equation}	
	Additionally, for every $i \neq j \in \Omega$, there exist scalars $\gamma_{i,j}$ for which
	\begin{equation}\label{condition2P}
	\gamma_{i,j} \bK_{j,i} \cP_j^{-1} \bK_{j,i}^T < \cP_i^{-1}.
	\end{equation}  
\end{assumption}

\begin{lemma}
	Given an LSS $\Si$ as defined in (\ref{LSS_def}), consider that the equations in (\ref{PgramLyap_gen}) have positive definite solutions $\cQ_q > \bfz, \ q \in \Omega$. Then, there exists $\mu>0$ such that for any switching signal in  (\ref{switch_signal}), with minimal dwell time $\mu$ (i.e. $t_i \geqslant \mu$) and $\bx(0) = \bfz$, the following bound holds 
	\begin{equation}\label{energy_bound_P}
	\bx^T(T_{\ell}) \cP_{q_\ell}^{-1} \bx(T_{\ell}) \leqslant \int_0^{T_\ell} \bu^T(s) \bu(s) ds.
	\end{equation}
\end{lemma}
\noindent
{\bf{Proof of Lemma 2}}. Consider that the conditions stated in Assumption 42 hold. Introduce $\gamma = \min\limits_{i,j \in \Omega, \ i \neq j}\gamma_{i,j}$ and let $\mu = - \frac{\ln \gamma}{M}$. For any piecewise continuous switching signal $\sigma: \mathbb{R} \rightarrow \Omega$ satisfying the conditions in (\ref{switch_signal}) and with minimal dwell time $\mu$, we will prove the bound stated in (\ref{energy_bound_P}). From (\ref{PgramLyap_gen}) and (\ref{condition1P}), it follows that
\begin{equation*}
\bA_i \cP_i +\cP_i \bA_i^T+ M_i \cP_i + \bB_i \bB_i^T \leqslant \bfz,
\end{equation*}
and by denoting $M =  \min\limits_{i \in \Omega}M_i$, the following holds
\begin{equation}\label{ineqP}
\bA_i \cP_i +\cP_i \bA_i^T+ M \cP_i + \bB_i \bB_i^T \leqslant \bfz.
\end{equation}
By multiplying the inequality (\ref{ineqP}) with $\cP_i^{-1}$ both to the left and to the right, we write
\begin{equation}\label{ineqP2}
\bA_i^T \cP_i^{-1} +\cP_i^{-1} \bA_i+ M \cP_i^{-1} + \cP_i^{-1} \bB_i \bB_i^T  \cP_i^{-1} \leqslant \bfz.
\end{equation}
Let $\bx(t)$ be the corresponding solution to (\ref{LSS_def}), and  also introduce the function  $V:\mathbb{R}^{n_i} \rightarrow \mathbb{R}$ as
\begin{equation}\label{Vdef2}
V(\bx(t)) = \begin{cases} \bx^T(t) \cP_{q_1}^{-1} \bx(t), \  t \in [0,t_1], \\ \bx^T(t) \cP_{q_i}^{-1} \bx(t),\  t \in (T_{i-1}, T_i], \ i \geqslant 2 \end{cases}.
\end{equation}
Using that $\dot{\bx}(t) = \bA_{q_i} \bx(t) + \bB_{q_i} \bu(t)$ and the definition of $V(\bx(t))$ in (\ref{Vdef2}), for $t \in (T_{i-1},T_i]$, we have
\begin{align*}
\frac{\partial V(\bx(t))}{\partial t}&= \frac{d \bx^T(t)}{d t} \cP_{q_i} ^{-1} \bx(t) +  \bx^T(t) \cP_{q_i}^{-1} \frac{d \bx(t)}{d t} = \bx^T(t) \big{(} \bA_{q_i}^T \cP_{q_i}^{-1} + \cP_{q_i}^{-1} \bA_{q_i}  \big{)} \bx(t) \\
&+ 2 \bx(t)^T \cP_{q_i}^{-1} \bB_{q_i} \bu(t),
\end{align*}
and by using the inequality in (\ref{ineqP2}), it follows that
\begin{align}
\frac{\partial V(\bx(t))}{\partial t}+ M V(\bx(t)) &\leqslant -\bx(t)^T \cP_{q_i}^{-1} \bB_{q_i} \bB_{q_i}^T \cP_{q_i}^{-1} \bx(t) + 2 \bx(t)^T \cP_{q_i}^{-1} \bB_{q_i} \bu(t) \nonumber \\ &= -\Vert \bB_{q_i}^T \cP_{q_i}^{-1} \bx(t) - \bu(t) \Vert_2^2 + \bu(t)^T \bu(t) .
\end{align}
Hence, the following inequality holds as,
\begin{equation}\label{Vineq}
\frac{\partial V(\bx(t))}{\partial t}+ M V(\bx(t)) \leqslant \bu(t)^T \bu(t), \ \ t \in (T_{i-1},T_i].
\end{equation}
By denoting $W(\bx(t)) = e^{M(t-T_{i})} V(\bx(t))$, \ for $t \in (T_{i-1},T_i] $, it follows that 
\begin{equation}\label{Wineq}
\frac{\partial W(\bx(t))}{\partial t} = e^{M(t-T_{i})} \Big{(}  \frac{\partial V(\bx(t))}{\partial t}+ M V(\bx(t)) \Big{)},
\end{equation}
and by combining (\ref{Vineq}) and (\ref{Wineq}) and integrating from $T_{i-1}$ to t, we obtain
\begin{equation}\label{Wdiff2}
W(\bx(t))-W(\bx(T_{i-1}^+)) \leqslant  \int\limits_{T_{i-1}}^t e^{M(s-T_{i})} \bu^T(s) \bu(s) ds .
\end{equation}
Following the same line of thought as in Section 4.1.1, one can show that $W(\bx(T_i^{+})) \leqslant W(\bx(T_i))$. By combining this statement with the inequality in (\ref{Wdiff2}), and by using the fact that $e^{M(s-T_{i})} \leqslant 1, \ \forall s \in (T_{i-1},T_i]$, one can write
\begin{align}\label{Wineq_cont1}
& W(\bx(T_i)) - W(\bx(T_{i-1})) \leqslant   \int\limits_{T_{i-1}}^{T_i} e^{M(s-T_{i})} \bu^T(s) \bu(s) ds \leqslant  \int\limits_{T_{i-1}}^{T_i} \bu^T(s) \bu(s) ds \nonumber \\ & \text{since} \ s-T_i \leqslant 0 
\Rightarrow \sum\limits_{i=1}^\ell W(\bx(T_i)) - W(\bx(T_{i-1})) \leqslant    \sum\limits_{i=1}^\ell  \int\limits_{T_{i-1}}^{T_i}  \bu^T(s) \bu(s) ds  \nonumber \\ &  \Rightarrow W(\bx(T_\ell)) - W(\bx(0)) \leqslant  \int\limits_{0}^{T_\ell}  \bu^T(s) \bu(s) ds.
\end{align}
Since $\bx(0) = \bfz$, it follows that $W(\bx(0)) = \bfz$. Also, from the definition of the function W, it is clear that $W(\bx(T_\ell)) = V(\bx(T_\ell)) = \bx^T(T_\ell) \cP_\ell^{-1} \bx(T_\ell)$. Hence, from (\ref{Wineq_cont1}), we directly conclude that
\begin{equation}\label{ineqP_Tl}
\bx^T(T_\ell) \cP_{q_\ell}^{-1} \bx(T_\ell) \leqslant \int\limits_{0}^{T_\ell}  \bu^T(s) \bu(s) ds, \ \forall \ell \geqslant 1,
\end{equation}
which proves the result in (\ref{energy_bound_P}). \sq

\subsection{Balancing transformation and truncation}

In this section, we introduce the procedure for model order
reduction by balanced truncation, and we prove a bound
of the approximation error.

\begin{procedure}
	{\rm
		Let $\Si = (n_1,n_2,\ldots,n_D, \{(\bA_q, \bB_q, \bC_q)| q \in \Omega \},\{ \bK_{q_{i},q_{i+1}} | q_i, q_{i+1} \in \Omega \}, \bx_0)$ be a linear switched system. A balanced realization of $\Si$ is denoted with the similar notation $\bar{\Si} = (n_1,n_2,\ldots,n_D, \{(\bar{\bA}_q, \bar{\bB}_q, \bar{\bC}_q)| q \in \Omega\},\{ \bar{\bK}_{q_{i},q_{i+1}} | q_i, q_{i+1} \in \Omega \}, \bx_0)$ and can be constructed as follows
		\begin{enumerate}
			\item Compute positive solutions $\cP_q >0$ for the equations in (\ref{PgramLyap_gen}) as well as solutions $\cQ_q >0$ for the equations in (\ref{QgramLyap_gen}), where $q \in \{1,2,\ldots,D\}$.
			\item Find square factor matrices $\bU_q$ so that $\cP_q = \bU_q \bU_q^T$. Additionally, compute the eigenvalue decomposition of the symmetric matrix $\bU_q^T \cQ_q \bU_q$, as 
			\begin{equation*}
			\bU_q^T \cQ_q \bU_q = \bV_q \bLa_q^2 \bV_q^T,
			\end{equation*}
			where $\bLa_q$ is a diagonal matrix with the real entries sorted in decreasing order.
			\item Construct the transformation matrices $\bS_q \in \mathbb{R}^{n_q \times n_q}$ as follows
			\begin{equation}\label{defSq}
			\bS_q = \bLa_q^{1/2} \bV_q^T \bU_q^{-1}.
			\end{equation}
			\item The matrices corresponding to the balanced realization $\bar{\Si}$ are computed in the following way (for any $q,q_1,q_2 \in \Omega$)
			\begin{equation}
			\bar{\bA}_q = \bS_q \bA_q \bS_q^{-1}, \ \ \bar{\bB}_q = \bS_q \bB_q, \ \ \bar{\bC}_q = \bC_q \bS_q^{-1}, \ \ \bar{\bK}_{q_1,q_2} = \bS_{q_2} \bK_{q_1,q_2} \bS_{q_1}^{-1}.
			\end{equation}
		\end{enumerate}
	}
\end{procedure}
The reachability and observability transformed Gramians $\bar{\cP}_q$ and respectively $\bar{\cQ}_q$, corresponding to mode q, are equal to each other and equal to $\bLa_q$
\begin{equation}
\bar{\cP}_q = \bS_q \cP_q \bS_q^T = \bLa_q, \ \ \bar{\cQ}_q = \big{(} \bS_q^{-1} \big{)}^T \cQ_q \bS_q^{-1} = \bLa_q,
\end{equation}
To prove these results, proceed as follows
\begin{equation*}
\bS_q \cP_q \bS_q^T = \big{(} \bLa_q^{1/2} \bV_q^T \bU_q^{-1} \big{)} \big{(} \bU_q \bU_q^T \big{)} \big{(}  \bLa_q^{1/2} \bV_q^T \bU_q^{-1} \big{)}^T = \bLa_q^{1/2} \bV_q^T \bV_q \bLa_q^{1/2} = \bLa_q,
\end{equation*}
and similarly for the observability transformed Gramian. The following result holds for any $i \in \Omega$:
\begin{align}
\bar{\bA}_i \bLa_i +\bLa_i \bar{\bA}_i^T+ \sum_{j=1, \ j \neq i}^D \bar{\bK}_{j,i} \bLa_j \bar{\bK}_{j,i}^T+ \bar{\bB}_i \bar{\bB}_i^T = \bfz, \label{Lyap_reach_bal} \\
\bar{\bA}_i^T \bLa_i +\bLa_i \bar{\bA}_i+ \sum_{j=1, \ j \neq i}^D \bar{\bK}_{i,j}^T \bLa_j \bar{\bK}_{i,j}+ \bar{\bC}_i^T \bar{\bC}_i = \bfz. \label{Lyap_obser_bal}
\end{align}
We will prove only the first equality since the proof for the second is similar. By multiplying the equation in (\ref{PgramLyap_gen}) corresponding to mode i with $\bS_i$ to the left and with $\bS_i^T$ to the right, we write
\small
\begin{align*}
& \bS_i \bA_i \cP_i \bS_i^T + \bS_i \cP_i \bA_i^T \bS_i^T+ \sum_{j=1, \ j \neq i}^D \bS_i \bK_{j,i} \cP_j \bK_{j,i}^T \bS_i^T+ \bS_i \bB_i \bB_i^T \bS_i^T = \bfz 
\\ & \Rightarrow \big{(} \bS_i \bA_i \bS_i^{-1} \big{)} \big{(} \bS_i  \cP_i \bS_i^T \big{)}  + \big{(} \bS_i \cP_i \bS_i^{-1} \big{)} \big{(} (\bS_i^{-1})^T \bA_i^T \bS_i^T \big{)} + \sum_{j=1, \ j \neq i}^D \big{(} \bS_i \bK_{j,i} \bS_j^{-1} \big{)} \big{(} \bS_j  \cP_j \bS_j^T \big{)} \\ & \big{(} (\bS_j^{-1})^T \bK_{j,i}^T \bS_i^T\big{)}+ \bS_i \bB_i \bB_i^T \bS_i^T = \bfz  \Rightarrow \bar{\bA}_i \bLa_i +\bLa_i \bar{\bA}_i^T+ \sum_{j=1, \ j \neq i}^D \bar{\bK}_{j,i} \bLa_j \bar{\bK}_{j,i}^T+ \bar{\bB}_i \bar{\bB}_i^T = \bfz.
\end{align*}
\normalsize
\noindent
After the system is rewritten in the equivalent balanced realization, the next step will be to construct a reduced order system by eliminating states similar as to the linear case with no switching.
One can partition the balanced realization of the original LSS $\Si$ in the following way
\begin{equation}\label{partition_bal}
\bar{\bA}_i = \left[ \begin{array}{cc} \bar{\bA}_i^{11} & \bar{\bA}_i^{12} \\ \bar{\bA}_i^{21} & \bar{\bA}_i^{22}
\end{array} \right], \ \ \bar{\bB}_i = \left[ \begin{array}{c} \bar{\bB}_i^{1}  \\ \bar{\bB}_i^{2}
\end{array} \right], \ \ \bar{\bC}_i = \left[ \begin{array}{cc} \bar{\bC}_i^{1} & \bar{\bC}_i^{2}
\end{array} \right], \ \ \bar{\bK}_{i,j} = \left[ \begin{array}{cc} \bar{\bK}_{i,j}^{11} & \bar{\bK}_{i,j}^{12} \\ \bar{\bK}_{i,j}^{21} & \bar{\bK}_{i,j}^{22}
\end{array} \right],
\end{equation}
where $\bar{\bA}_i^{11} \in \mathbb{R}^{r_i \times r_i}, \ \bar{\bK}_{i,j}^{11} \mathbb{R}^{r_j \times r_i},  \ \bar{\bB}_i^{1} \in \mathbb{R}^{r_i}, \ \bar{\bC}_i^{1} \in \mathbb{R}^{1 \times r_i}$. The truncation orders are chosen to be less than the dimensions of the subsystems, i.e. $r_i \leqslant n_i$.
\begin{definition}
	Consider as given an original LSS $\Si$ and the balanced equivalent system $\bar{\Si}$ corresponding to $\Si$ for which the system matrices are split as in (\ref{partition_bal}). Let $\hat{\Si} = (r_1,r_2,\ldots,r_D, \{(\hat{\bA}_q, \hat{\bB}_q, \hat{\bC}_q)| q \in \Omega\},\{ \hat{\bK}_{q_{i},q_{i+1}} | q_i, q_{i+1} \in \Omega \}, \bx_0)$, be a reduced linear switched system for which the system matrices are written as follows
	\begin{equation}\label{reduced_bal}
	\hat{\bA}_q = \bar{\bA}_q^{11}, \ \ \hat{\bB}_q = \bar{\bB}_i^{1}, \ \ \hat{\bC}_q = \bar{\bC}_q^{1}, \ \ \hat{\bK}_{q_1,q_2} = \bar{\bK}_{q_1,q_2}^{11},
	\end{equation} 
	where $r_q \leqslant n_q$ and $q,q_1,q_2 \in \Omega$.
\end{definition}
By writing the dynamics of both the original balanced system $\bar{\Si}$ and the reduced system $\hat{\Si}$, as
\begin{equation}\label{dynamics_bal}
\dot{\bar{\bx}}(t) = \bar{\bA}_{q_i} \bar{\bx}(t) + \bar{\bB}_{q_i} \bu(t), \ \ \dot{\hat{\bx}}(t) = \hat{\bA}_{q_i} \hat{\bx}(t) + \hat{\bB}_{q_i} \bu(t), \ \ t \in (T_{i-1},T_i],
\end{equation}
and continuing with the transition of the state variable from mode $q_i$ to mode $q_{i+1}$ at time $T_i$ again for both systems
\begin{equation}\label{bxtildebar}
\bar{\bx}(T_i^{+}) = \bar{\bK}_{q_i,q_{i+1}} \bar{\bx}(T_i), \ \ \ \hat{\bx}(T_i^{+}) = \hat{\bK}_{q_i,q_{i+1}} \hat{\bx}(T_i),
\end{equation}
we finally conclude that the original output and the one corresponding to the reduced LSS are written as
\begin{equation}
\bar{\by}(t) = \bar{\bC}_{q_i} \bar{\bx}(t) = \bC_{q_i} \bx(t) = \by(t), \ \ \hat{\by}(t) = \hat{\bC}_{q_i} \hat{\bx}(t).
\end{equation}
We also partition the balanced Gramians corresponding to the system $\bar{\Si}$ as
\begin{equation}\label{gramm_split}
\bLa_i = \left[ \begin{array}{cc} \hat{\bLa}_i & 0 \\ 0 & \check{\bLa}_i
\end{array} \right],\ \  \hat{\bLa}_i \in \mathbb{R}^{r_i}, \ \check{\bLa}_i \in \mathbb{R}^{n_i-r_i}.
\end{equation}
By plugging in the matrices in (\ref{partition_bal}) and (\ref{gramm_split}), into the equation (\ref{Lyap_reach_bal}), it follows that
\begin{equation}
\hat{\bA}_i \hat{\bLa}_i +\hat{\bLa}_i \hat{\bA}_i^T+ \sum_{j=1, \ j \neq i}^D \hat{\bK}_{j,i} \hat{\bLa}_j \hat{\bK}_{j,i}^T+ \sum_{j=1, \ j \neq i}^D \bar{\bK}^{12}_{j,i} \check{\bLa}_j \big{(} \bar{\bK}^{12}_{j,i} \big{)}^T+ \hat{\bB}_i \hat{\bB}_i^T = \bfz_{r_i}.
\end{equation}
Note that the reduced balanced Gramians $\hat{\bLa}_i$ do not satisfy the same type of generalized Lyapunov equations as the original balanced Gramians, i.e., the equations in (\ref{Lyap_reach_bal}) and  (\ref{Lyap_obser_bal}). Instead,  conclude that the diagonal Gramians $\hat{\bLa}_i$ satisfy the following inequalities 
\begin{align}
\hat{\bA}_i \hat{\bLa}_i +\hat{\bLa}_i \hat{\bA}_i^T+ \sum_{j=1, \ j \neq i}^D \hat{\bK}_{j,i} \hat{\bLa}_j \hat{\bK}_{j,i}^T+ \hat{\bB}_i \hat{\bB}_i^T < \bfz_{r_i} \\ 
\hat{\bA}_i^T \hat{\bLa}_i +\hat{\bLa}_i \hat{\bA}_i+ \sum_{j=1, \ j \neq i}^D \hat{\bK}_{i,j}^T \hat{\bLa}_j \hat{\bK}_{i,j}+ \hat{\bC}_i^T \hat{\bC}_i < \bfz_{r_i}.
\end{align}
Hence, the reduced-order diagonal matrices $\hat{\bLa}_i, \ i \in \Omega$ could be also considered as Gramians, but in the relaxed way introduced in Definition 12.

\subsubsection{Error bound}

In this section we present a bound on the $L_2$ norm of the difference between the observed outputs corresponding to the original LSS and to the reduced LSS. We will show that this bound depends on the $L_2$ norm of the chosen control input and on the neglected elements of the balanced reduced Gramians. Some of the derivations presented here are inspired from techniques used prior in the dissertations \cite{sandberg_phd,besselink_phd} and in the more recent contribution \cite{bdc17}, that provides a bound for BT applied to stochastic systems. 

We assume that all pairs of the original Gramians $(\cP_i,\cQ_i)$, defined as the solutions of the equations (\ref{PgramLyap_gen}) and (\ref{QgramLyap_gen}), are transformed through the corresponding balanced transformations $\bV_i$, into $(\bLa_i,\bLa_i)$ where $\bLa_i$ are diagonal matrices ($i \in \Omega$).

\begin{assumption} 
	Consider that $\sum_{j=1, \ j \neq i}^D \bar{\bK}_{j,i} \bLa_j \bar{\bK}_{j,i}^T >0$ and $\sum_{j=1, \ j \neq i}^D  \bar{\bK}_{i,j}^T \bLa_j$ $\bar{\bK}_{i,j} > 0$. Hence, one can always choose an $M > 0$ such that the following conditions hold
	\begin{equation} \label{condition1L}
	\sum_{j=1, \ j \neq i}^{D} \bar{\bK}_{j,i} \bLa_j \bar{\bK}_{j,i}^T > M \bLa_i, \ \sum_{j=1, \ j \neq i}^{D} \bar{\bK}_{i,j}^T \bLa_j \bar{\bK}_{i,j} > M \bLa_i, \ \forall i \in \Omega.
	\end{equation}	
\end{assumption} 
\noindent
By replacing inequalities in (\ref{condition1L}), into equations (\ref{Lyap_reach_bal}) and (\ref{Lyap_obser_bal}), it follows that
\begin{equation}
\bar{\bA}_i \bLa_i +\bLa_i \bar{\bA}_i^T+ M \bLa_i+ \bar{\bB}_i \bar{\bB}_i^T < \bfz, \ \ 
\bar{\bA}_i^T \bLa_i +\bLa_i \bar{\bA}_i+ M \bLa_i+ \bar{\bC}_i^T \bar{\bC}_i < \bfz. \label{Lyap_ineq_bal}
\end{equation}
By multiplying the first inequality in (\ref{Lyap_ineq_bal}) with $\bLa^{-1}$ to the left and to the right, one can again write that
\begin{equation}\label{Lyap_ineq_bal_reach}
\bar{\bA}_i^T \bLa_i^{-1} +\bLa_i^{-1} \bar{\bA}_i+ M \bLa_i^{-1}+ \bLa_i^{-1} \bar{\bB}_i \bar{\bB}_i^T \bLa_i^{-1} < \bfz.
\end{equation}
From (\ref{Lyap_ineq_bal_reach}) and the second inequality in (\ref{Lyap_ineq_bal}), it directly follows that the following relations hold for any vectors $\bz$ and $\bv$
\begin{align}
2 (\bar{\bA}_i \bz+ \bar{\bB}_i \bv ) \bLa_i^{-1} \bx  \leqslant  \Vert \bv \Vert_2^2 - M \bz^T \bLa_i^{-1} \bz \label{ineqv} , \\
2 \bz^T \bar{\bA}_i^T \bLa_i \bz \leqslant - \Vert \bar{\bC}_i \bz \Vert_2^2 - M \bz^T \bLa_i \bz. \label{ineqz}
\end{align}
Next, for all $i \in \{1,2,\ldots,D\}$, proceed to partition the transformed Gramians $\bLa_i$ 
\begin{equation}
\bLa_i = \left[ \begin{array}{cc} \hat{\bLa}_i & 0 \\ 0 & \beta_i
\end{array} \right],\ \  \beta_i \in \mathbb{R}.
\end{equation}
Let $\beta = \max\limits_{i \in \Omega}{\beta_i}$. By splitting the state variable $\bar{\bx}(t)$ as $\bar{\bx}(t) = \left[ \begin{array}{cc} \bar{\bx}_1(t) & \bar{\bx}_2(t) \end{array} \right]^T$, with $\bar{\bx}_1(t) \in \mathbb{R}^{n-1}, \ \bar{\bx}_2(t) \in \mathbb{R}$, introduce real valued vectors
\begin{equation}\label{xoxc}
\bx_o(t) = \left[ \begin{array}{c} \bar{\bx}_1(t) - \hat{\bx}(t) \\ \bar{\bx}_2(t) \end{array} \right], \ \ \ \bx_c(t) = \left[ \begin{array}{c} \bar{\bx}_1(t) + \hat{\bx}(t) \\ \bar{\bx}_2(t) \end{array} \right].
\end{equation}
Note that the following holds:
\begin{equation*}
\by(t) - \hat{\by}(t) = \bC_{q_i} \bx_o(t), \ \ t \in (T_{i-1}, T_i].
\end{equation*}
Define the function $V: \mathbb{R}^n \times \mathbb{R}^n \rightarrow \mathbb{R}$ as follows
\begin{equation}\label{defVxoxc}
V(\bx_o(t), \bx_c(t)) = \bx_o(t)^T \bLa_{q_i} \bx_o(t)+ \beta_{q_i}^2 \bx_c(t)^T(t) \bLa_{q_i}^{-1} \bx_c(t), \  t \in (T_{i-1}, T_i].
\end{equation}

\begin{lemma} The temporal derivative of the function V, as defined in (\ref{defVxoxc}), satisfies
	\begin{equation}\label{derVineq}
	\frac{\partial V(\bx_o(t), \bx_c(t))}{\partial t} \leqslant - M V(t) + 4 \beta^2 \Vert \bu(t) \Vert_2^2- \Vert \by(t)-\hat{\by}(t) \Vert_2^2.
	\end{equation}
\end{lemma}
{\bf{Proof of Lemma 3}}. By putting together (\ref{partition_bal}), (\ref{reduced_bal}) and (\ref{dynamics_bal}) and by using the notation in (\ref{xoxc}), we can write that
\begin{align}
\dot{\bx}_o(t) &= \bA_{q_i} \bx_o (t) + \left[ \begin{array}{c} \bfz \\ \bB_{q_i}^2(t) \end{array} \right] \bu(t) + \left[ \begin{array}{c} \bfz \\ \bA_{q_i}^{21}(t) \end{array} \right] \bar{\bx}(t), \label{xo_der} \\
\dot{\bx}_c(t) &= \bA_{q_i} \bx_c (t) +2 \bB_{q_i}^2 \bu(t) - \left[ \begin{array}{c} \bfz \\ \bB_{q_i}^2(t) \end{array} \right] \bu(t) + \left[ \begin{array}{c} \bfz \\ \bA_{q_i}^{21}(t) \end{array} \right] \bar{\bx}(t). \label{xc_der}
\end{align}
By using (\ref{xo_der}) and the inequality in (\ref{ineqz}), one can write that
\small
\begin{align}\label{derVxo}
& \frac{d }{d t} \bx_o(t)^T \bLa_{q_i} \bx_o(t) = 2 \bx_o^T (t) \bLa_{q_i} \bx_o(t) + 2 \Big{(} \left[ \begin{array}{c} \bfz \\ \bB_{q_i}^2 \bu(t)+\bA_{q_i}^{21} \hat{\bx}(t) \end{array} \right]^T \bLa_{q_i} \bx_o(t) \Big{)} \nonumber \\
& \leqslant -M \bx_o^T(t) \bLa_{q_i} \bx_o(t)- \Vert \bC_{q_i} \bx_o(t)  \Vert_2^2 + 2 \alpha_o  = -M \bx_o^T(t) \bLa_{q_i}^{-1} \bx_o(t) - \Vert \by(t) - \hat{\by}(t) \Vert_2^2  +  2 \alpha_o,
\end{align}
\normalsize
where  
\small
\begin{align} \label{alphao}
\alpha_o  = \left[ \begin{array}{c} \bfz \\ \bB_{q_i}^2 \bu(t)+\bA_{q_i}^{21} \hat{\bx}(t) \end{array} \right]^T  \left[ \begin{array}{cc} \hat{\bLa}_{q_i} & \bfz \\ \bfz & \beta_{q_i} \end{array} \right] \left[ \begin{array}{c} \bar{\bx}_1(t)- \hat{\bx}(t) \\ \bar{\bx}_2(t) \end{array} \right]  =  \beta_{q_i} \big{(} \bB_{q_i}^2 \bu(t)+\bA_{q_i}^{21} \hat{\bx}(t) \big{)}^T \bar{\bx}_2(t)
\end{align}
\normalsize
Similarly, by using (\ref{xc_der}) and the inequality in (\ref{ineqv}) for $\bz = \bx_c(t)$ and $\bv = 2 \bu(t)$ , one can show that
\small
\begin{align}\label{derVxc}
\frac{d }{d t} \bx_c(t)^T \bLa_{q_i}^{-1} \bx_c(t) &= 2 \big{(} \bA_{q_i}  \bx_c(t) + \bB_{q_i} 2 \bu(t)  \big{)} \bLa_{q_i}^{-1} \bx_o(t) - 2 \Big{(} \left[ \begin{array}{c} \bfz \\ \bB_{q_i}^2 \bu(t)+\bar{\bA}_{q_i}^{21} \hat{\bx}(t) \end{array} \right]^T \bLa_{q_i}^{-1} \bx_c(t) \Big{)} \nonumber \\
& \leqslant -M \bx_c^T(t) \bLa_{q_i}^{-1} \bx_c(t)+ 4 \Vert \bu(t) \Vert_2^2- 2 \alpha_c ,
\end{align}
\normalsize
where 
\small
\begin{align} \label{alphac}
\alpha_c  = \left[ \begin{array}{c} \bfz \\ \bB_{q_i}^2 \bu(t)+\bA_{q_i}^{21} \bar{\bx}(t) \end{array} \right]^T  \left[ \begin{array}{cc} \hat{\bLa}_{q_i}^{-1} & \bfz \\ \bfz & \beta_{q_i}^{-1}  \end{array} \right] \left[ \begin{array}{c} \bar{\bx}_1(t)+ \hat{\bx}(t) \\ \bar{\bx}_2(t) \end{array} \right]  =  \beta_{q_i}^{-1} \big{(} \bB_{q_i}^2 \bu(t)+\bA_{q_i}^{21} \bar{\bx}(t) \big{)}^T \hat{\bx}_2(t)
\end{align}
\normalsize
From (\ref{alphao}) and (\ref{alphac}), observe that $\alpha_o = \beta_{q_i}^2 \alpha_c$. Hence, by adding the inequality in (\ref{derVxo}) with the one in (\ref{derVxc}) multiplied by $\beta_{q_i}^2$, it follows that
\begin{align*}
& \frac{d }{d t} \bx_o(t)^T \bLa_{q_i} \bx_o(t) + \beta_{q_i}^2  \frac{d }{d t} \bx_c(t)^T \bLa_{q_i}^{-1} \bx_c(t) \leqslant -M \big{(}   \bx_o(t)^T \bLa_{q_i} \bx_o(t) \\ &+ \beta_{q_i}^2  \bx_c(t)^T \bLa_{q_i}^{-1} \bx_c(t) \big{)}   - \Vert \by(t) - \hat{\by}(t) \Vert_2^2 + 4 \beta_{q_i}^2 \Vert \bu(t) \Vert_2^2,
\end{align*}
and by using the definition of $V(t)$ in (\ref{defVxoxc}), it automatically proves the result in (\ref{derVineq}). $\sq$

Introduce the concatenation of the state variables and of the coupling matrices corresponding to the (balanced) original and reduced systems,  $t \in (T_{\ell-1},T_\ell]$
\begin{equation}\label{defxhat}
\tilde{\bx}(t) = \left[ \begin{array}{c} \bar{\bx}(t) \\   \hat{\bx}(t)  \end{array} \right] \in \mathbb{R}^{2n_{q_i}-1}, \ \
\tilde{\bK}_{q_i,q_{i+1}} = \left[ \begin{array}{cc} \bar{\bK}_{q_i,q_{i+1}} & \bfz \\  \bfz & \hat{\bK}_{q_i,q_{i+1}}  \end{array} \right] \in \mathbb{R}^{2n_{q_i}-1 \times 2n_{q_i}-1}. 
\end{equation}
From (\ref{bxtildebar}) and (\ref{defxhat}), it follows that $\tilde{\bx}(T_i^+) = \tilde{\bK}_{q_i,q_{i+1}}   \tilde{\bx}(T_i)$. Note that the function $V: \mathbb{R}^n \times \mathbb{R}^n \rightarrow \mathbb{R}$ defined in (\ref{defVxoxc}), can also be written as a function of $\tilde{\bx}(t)$, as
\begin{equation}\label{defVxhat}
V(\tilde{\bx}(t)) = \tilde{\bx}(t)^T \tilde{\bR}_{q_i} \tilde{\bx}(t) = \left[ \begin{array}{c} \bar{\bx}(t) \\  \hat{\bx}(t) \end{array} \right]^T \tilde{\bR}_{q_i} \left[ \begin{array}{c} \bar{\bx}(t) \\  \hat{\bx}(t) \end{array} \right], \ \ t \in (T_{i-1},T_i],
\end{equation}
where the matrices $\tilde{\bR}_q \in \mathbb{R}^{2n_q-1 \times 2n_q-1}$ are defined for any $q \in \Omega$, as
\small
\begin{align}\label{defRq}
\tilde{\bR}_q = \left[ \begin{array}{ccc} \hat{\bLa}_q & \bfz & -\hat{\bLa}_q \\ 0 & \beta_q & 0 \\ -\hat{\bLa}_q & \bfz & \hat{\bLa}_q  \end{array} \right]+ \beta_q^2 \left[ \begin{array}{ccc} \hat{\bLa}_q^{-1} & \bfz & \hat{\bLa}_q^{-1} \\ 0 & \beta_q & 0 \\ \hat{\bLa}_q^{-1} & \bfz & \hat{\bLa}_q^{-1} \end{array} \right] = \left[ \begin{array}{ccc} \hat{\bLa}_q+\beta_q^2 \hat{\bLa}_q^{-1} & \bfz & -\hat{\bLa}_q+\beta_q^2 \hat{\bLa}_q^{-1} \\ 0 & 2\beta_q & 0 \\ -\hat{\bLa}_q+\beta_q^2 \hat{\bLa}_q^{-1} & \bfz & \hat{\bLa}_q+\beta_q^2 \hat{\bLa}_q^{-1} \end{array} \right].
\end{align}
\normalsize
\begin{assumption}
	Since $\tilde{\bR}_q > \bfz$, one can consider that there exists $\gamma \in (0,1)$ so that for all $q_i, q_{i+1} \in \Omega$, so that
	\begin{equation}
	\gamma \tilde{\bK}_{q_i,q_{i+1}}^T \tilde{\bR}_{q_{i+1}} \tilde{\bK}_{q_i,q_{i+1}} < \tilde{\bR}_{q_{i}}.
	\end{equation}
\end{assumption}
First, we present a result for one step reduction. The $L_2$ norm of the output error computed as the differences between the original output and the output corresponding to the reduced system is bounded by the norm of the input.
\begin{theorem}
	Let $\Si = (n_1,n_2,\ldots,n_D, \{(\bA_q, \bB_q, \bC_q)| q \in \Omega \},\{ \bK_{q_{i},q_{i+1}} | q_i, q_{i+1} \in \Omega \}, \bx_0)$ be a linear switched system and let $\hat{\Si}$ be a reduced order system obtained from $\Si$ via the proposed balancing and truncation procedure,
	$$
	\hat{\Si} = (n_1-1,n_2-1,\ldots,n_D-1, \{(\hat{\bA}_q, \hat{\bB}_q, \hat{\bC}_q)| q \in \Omega\},\{ \hat{\bK}_{q_{i},q_{i+1}} | q_i, q_{i+1} \in \Omega \}, \hat{\bx}_0).
	$$
	Consider any control inputs $\bu(t) \in L_2(\mathcal{R}^m)$ and denote with $\by(t)$ and $\hat{\by}(t)$ the outputs of the systems $\Si$ and, respectively $\hat{\Si}$ for the zero state case (i.e., $\bx(0) = \bfz$). Then, there exists $\mu >0$ such that for any switching signal with minimal dwell time $\mu$ (i.e. $t_i > \mu, \ \forall i$), so that
	\begin{equation}\label{errbound}
	\Vert \by - \hat{\by} \Vert_{2} \leqslant 2 \beta \Vert \bu \Vert_{2}.
	\end{equation}
\end{theorem}
\noindent
{\bf{Proof of Theorem 1}}. 
Choose $\mu = -\frac{\ln \gamma}{M}$ as the minimal dwell time for the switching signal $\sigma(t)$ (where $\gamma$ is as in Assumption 44).

Introduce the function $W(\tilde{\bx}(t)) = e^{M(t-T_{i-1})} V(\tilde{\bx}(t))$, \ for $t \in (T_{i-1},T_i] $. It follows that 
\begin{equation}\label{Wdef3}
\frac{\partial W(\tilde{\bx}(t))}{\partial t} = e^{M(t-T_{i-1})} \Big{(}  \frac{\partial V(\tilde{\bx}(t))}{\partial t}+ M V(\tilde{\bx}(t)) \Big{)}.
\end{equation}
Let $\Theta(t) = 4 \beta^2 \Vert \bu(t) \Vert_2^2- \Vert \by(t)-\bar{\by}(t) \Vert_2^2$. From (\ref{derVineq}) and (\ref{Wdef3}), we write that
\begin{equation}
\frac{\partial W(\tilde{\bx}(t))}{\partial t} \leqslant e^{M(t-T_{i-1})} \Theta(t), \ \ t \in (T_{i-1},T_i].
\end{equation}
Repeating the reasoning from the proof of Lemma 2, it follows that
\begin{equation}
\tilde{\bx}^T(T_\ell) \tilde{\bR}_{q_\ell} \tilde{\bx}(T_\ell) \leqslant \int_{0}^{T_\ell} \Theta(s) ds, \ \ \forall \ell \geqslant 1.
\end{equation}
Since $\hat{\bR}_{q_\ell} > 0$, then $\int_{0}^{T_\ell} \Theta(s) ds \geqslant 0, \ \ \forall \ell \geqslant 1$. 
By allowing $T_\ell \rightarrow \infty$ and by using the definition of the function $\Theta$, we can write
\begin{equation*}
4 \beta^2 \int_{0}^{\infty} \Vert \bu(s) \Vert_2^2 ds \geqslant \int_{0}^{\infty} \Vert \by(s)-\hat{\by}(s) \Vert_2^2 ds.
\end{equation*}
Hence, the result in (\ref{errbound}) has been proven. \sq

\vspace{2mm}

\noindent

By partitioning the set of discrete modes in two disjoint subsets, as $\Omega = \{1,2,\ldots,D\}$  $=  \Omega_1 \bigcup \Omega_2$, we emphasize two different cases when reducing the  system $\Si$
\begin{equation}\label{setup_new}
\begin{cases} q \in \Omega_1 \ \ \Rightarrow \ \ \text{perform reduction by 1 of the LTI subsystem in mode q}, \\ q \in \Omega_2 \ \ \Rightarrow \ \ \text{do not change the LTI in mode q}.  \end{cases}
\end{equation} 
Next, introduce the balanced Gramians corresponding to the two subsets, as
\begin{equation}
\bLa_\ell = \left[ \begin{array}{cc} \hat{\bLa}_\ell & \bfz \\ \bfz & \beta_\ell \end{array} \right], \ \text{for} \ \ell \in \Omega_1, \ \ \ \text{and} \ \ \ \hat{\bLa}_\ell = \bLa_{\ell}, \ \text{for} \ \ell \in \Omega_2.
\end{equation}
\begin{remark}
	We conclude that the bound in (\ref{errbound}) still holds for the setup that was introduced in (\ref{setup_new}), as follows 
	\begin{equation}\label{errbound_extended}
	\Vert \by - \hat{\by} \Vert_{2} \leqslant 2 \beta \Vert \bu \Vert_{2}, \ \ \ \beta = \max\limits_{\ell \in \Omega_1} \beta_\ell.
	\end{equation}
	Here, the selection of the scalar $\beta$ is restricted only to  diagonal Gramians  corresponding to the discrete modes from $\Omega_1$. The proof is similar to the one just presented and will be skipped for brevity reasons.
\end{remark}

Next, we will present a more general result by extending Theorem 1 from one step reduction to reduction to any dimension by allowing possibly different reduction levels for each active mode $q \in \Omega$. Consider that the diagonal Gramians associated to the original and reduced systems can be written as
\begin{equation}
\bLa_q = \left[ \begin{array}{ccc} \sigma_{q,1} &  & 0 \\  & \ddots &  \\ 0 &  & \sigma_{q,n_q}
\end{array} \right] \in \mathbb{R}^{n_q \times n_q}, \ \ \ \hat{\bLa}_q = \left[ \begin{array}{ccc} \sigma_{q,1} &  & 0 \\  & \ddots &  \\ 0 & & \sigma_{q,r_q} 
\end{array} \right] \in \mathbb{R}^{r_q \times r_q}.
\end{equation}
for $q \in \{1,2,\ldots,D\}$ and $\sigma_{q,1} \geqslant \sigma_{q,2} \geqslant \ldots \geqslant \sigma_{q,r_q} \geqslant \ldots \sigma_{q,n_q}>0 $.
For $\ell \in \{1,2,\ldots,\xi\}$, introduce the following diagonal matrices
\begin{equation}
_\ell\hat{\bLa}_q =  \left[ \begin{array}{ccc} \sigma_{q,1} &  & 0 \\  & \ddots &  \\ 0 &  & \sigma_{q,n_q-i+1} 
\end{array} \right], \ \text{if} \ i \leqslant n_q-r_q,  \   _\ell\hat{\bLa}_q =  \left[ \begin{array}{ccc} \sigma_{q,1} &  & 0 \\  & \ddots &  \\ 0 &  & \sigma_{q,r_q} 
\end{array} \right], \ \text{if} \ i> n_q-r_q 
\end{equation}
Similarly, let $_\ell\hat{\bA}_q  \in \mathbb{R}^{r_{q,\ell} \times r_{q,\ell}}, \ _\ell\hat{\bB}_q  \in \mathbb{R}^{r_{q,\ell} \times m_q}, \ _\ell\hat{\bC}_q \in  \mathbb{R}^{p_q \times r_{q,\ell} }, \ _\ell\hat{\bK}_{q_1,q_2} \in \mathbb{R}^{r_{q_2,\ell} \times r_{q_1,\ell}}$, be the $(1,1)$ blocks of the matrices defined in (\ref{partition_bal});
\small
\begin{equation}\label{def_lAq}
\bar{\bA}_q = \left[ \begin{tabular}{ l |  r }			
$_\ell\hat{\bA}_q$ & \ \ \  \  \\
\hline
&  \\
\end{tabular}  \right], \ \  \bar{\bB}_q = \left[ \begin{tabular}{ l   }			
$_\ell\hat{\bB}_q$  \\
\hline
\\
\end{tabular} \right], \ \ \bar{\bC}_q = \left[ \begin{tabular}{ l | r }			
$_\ell\hat{\bC}_q$ & \\
\end{tabular} \right], \ \ \bar{\bK}_{q_1,q_2} = \left[ \begin{tabular}{ l |  r }			
$_\ell\hat{\bK}_{q_1,q_2}$ & \ \ \ \ \ \  \\
\hline
&  \\
\end{tabular}  \right].
\end{equation}
\normalsize
for $r_{q,\ell} = \begin{cases} n_q - \ell, \ \ \text{if} \ \ell \leqslant n_q - r_q \\ \ \ \ \ r_q,  \ ~ \ \ \text{if} \ \ell > n_q - r_q   \end{cases}$.

\begin{definition}
	Using the matrices introduced in (\ref{def_lAq}), construct the family of reduced linear switched systems $\{\hat{\Si}_\ell \ \vert \ 0 \leqslant \ell \leqslant \xi \}$ with $\xi = \max\limits_{q \in \Omega}(n_q-r_q)$, as 
	\begin{equation}
	\hat{\Si}_{\ell} = (r_{1,\ell},r_{2,\ell},\ldots,r_{D,\ell}, \{(_\ell\hat{\bA}_q, _\ell\hat{\bB}_q, _\ell\hat{\bC}_q)| q \in \Omega\},\{ _\ell\hat{\bK}_{q_{i},q_{i+1}} | q_i, q_{i+1} \in \Omega \}, \hat{\bx}_0).
	\end{equation}
\end{definition}
\vspace{-2mm}
\noindent
Note that for $\ell = 0$, the element $\bar{\Si}_0$ coincides to the original LSS in balanced format, i.e. $\hat{\Si}_0 = \bar{\Si}$. Moreover, when $\ell = \xi$, it follows that $\hat{\Si}_\xi = \hat{\Si}$, with $\hat{\Si}$ as  introduced in Definition 13.

\begin{theorem}
	Let $\Si = (n_1,n_2,\ldots,n_D, \{(\bA_q, \bB_q, \bC_q)| q \in \Omega \},\{ \bK_{q_{i},q_{i+1}} | q_i, q_{i+1} \in \Omega \}, \bx_0)$ be a linear switched system and let $\hat{\Si}_\ell$ be a reduced order system obtained from $\Si$ introduced in Definition 14. Consider the control input $\bu(t) \in L_2(\mathbb{R}^m)$ and denote with $\by(t)$ and $\hat{\by}(t)$ the outputs of the systems $\Si$ and, respectively $\hat{\Si}$ for the zero state case (i.e., $\bx(0) = \bfz$). Then, there exists $\mu >0$ such that for any switching signal with minimal dwell time $\mu$ (i.e. $t_i > \mu, \ \forall i$), so that
	\begin{equation}\label{errbound2}
	\Vert \by - \hat{\by} \Vert_{2} \leqslant 2 \beta \Vert \bu \Vert_{2},
	\end{equation}
	where $\beta = \sum\limits_{\ell=1}^{\xi} \eta_\ell, \ \ \eta_\ell =  \max\limits_{\ell  \leqslant n_q-r_q, \ q \in \Omega} \sigma_{q,n_q-\ell+1}$.
\end{theorem}
\noindent
{\bf{Proof of Theorem 2}}. We start by applying the result of Theorem 1 (for one step reduction) as adapted in Remark 6 (allowing adjustable reduction levels for different modes), to $\hat{\Si}_{\ell-1}$ and $\hat{\Si}_{\ell}$ for all $\ell \in \{1,2,\ldots,\xi \}$. Consider the following two subsets of $\Omega$,
\begin{equation*}
\Omega_1^{\ell} = \{ q \in \Omega \ \vert \ell \leqslant n_q - r_q \}
,  \ \ \Omega_2^{\ell} = \{ q \in \Omega \ \vert \ell >  n_q - r_q \}.
\end{equation*}
Note that $\hat{\Si}_\ell$ is the result of a one-step reduction applied to $\hat{\Si}_{\ell-1}$. Also, the following inequalities hold for all $\ell \in \{1,2,\ldots,\xi\}$
\begin{align*}
_{\ell-1}\hat{\bA}_q  \ {_{\ell-1}}\hat{\bLa}_q + {_{\ell-1}} \hat{\bLa}_q \ {_{\ell-1}} \hat{\bA}_q^T+ M \ {_{\ell-1}} \hat{\bLa}_q+ {_{\ell-1}}\hat{\bB}_q \ {_{\ell-1}} \hat{\bB}_q^T < \bfz, \\ 
{_{\ell-1}}\hat{\bA}_q^T \ {_{\ell-1}} \hat{\bLa}_q + {_{\ell-1}} \hat{\bLa}_q \ {_{\ell-1}}\hat{\bA}_q+ M \ {_{\ell-1}} \hat{\bLa}_q+ {_{\ell-1}}\hat{\bC}_q^T \ {_{\ell-1}}\hat{\bC}_q < \bfz.
\end{align*}
where M is chosen as in the proof of Theorem 1 and it does not depend on $\ell$. Next, denote with $\hat{\by}_{\ell}$ and $\hat{\by}_{\ell-1}$ the outputs corresponding to the systems $\hat{\Si}_{\ell}$ and, respectively $\hat{\Si}_{\ell-1}$ for input $\bu(t) \in L_2$, switching signal $\sigma(t)$ with minimal dwell time $\mu_\ell$ and initial zero states. From (\ref{errbound_extended}), it follows that 
\begin{equation}\label{errbound3}
\Vert \hat{\by}_{\ell-1} - \hat{\by}_\ell \Vert_{2} \leqslant 2 \eta_{\ell} \Vert \bu \Vert_{2}.
\end{equation}
For $\ell = 0$, the output $\hat{\by}_0$ coincides to the output of the original LSS in balanced format, i.e. $\hat{\by}_0 = \by$. Furthermore, when $\ell = \xi$, it follows that $\hat{\by}_\xi = \hat{\by}$, with $\hat{\by}$ as   in Section 4.3, i.e. the output of the reduced-order LSS $\hat{\Si}$ from Definition 13. By adding the inequalities in (\ref{errbound3}) for all values of $\ell$ in $\{1,\ldots,\xi\}$, it follows that
\begin{align*} \label{errbound3}
\sum_{\ell = 1}^\xi \Vert \hat{\by}_{\ell-1} - \hat{\by}_\ell \Vert_{2} \leqslant 2 \sum_{\ell = 1}^\xi \eta_{\ell} \Vert \bu \Vert_{2} \Rightarrow \Vert \sum_{\ell = 1}^\xi  ( \hat{\by}_{\ell-1} - \hat{\by}_\ell) \Vert_{2} \leqslant 2 \beta \Vert \bu \Vert_{2} \\  \Rightarrow  \Vert \hat{\by}_{0} - \hat{\by}_\xi \Vert_{2} \leqslant 2 \beta \Vert \bu \Vert_{2}, \hspace{30mm}
\end{align*}
which implies that the result in (\ref{errbound2}) is thus proven. \sq

\begin{example}
	\rm {
		To clarify the notation used in the proof of Theorem 2, we present a simple example for $D = 3$, i.e. $\Omega = \{1,2,3\}$. Assume $n_q = 3, \ \forall q \in \Omega$ and the choose reduction orders 1,3 and respectively, 2 for modes 1,2 and respectively, 3. Also, note that $\xi = \max\limits_{q \in \Omega}(n_q-r_q) = 2$.
		
		\begin{align*}
		\begin{tabular}{|l|c|r|l|}
		\hline 
		\ \ \ \ $q$ & 1 & 2 & 3 \\
		\hline 
		\ \ \ $n_q$ & 3 & 3 & 3 \\
		\hline
		\ \ \ $r_q$ & 1 & 3 & 2 \\
		\hline
		$n_q-r_q$ & 2 & 0 & 1 \\
		\hline
		\end{tabular} \ \ \ \Rightarrow \ \ \ \begin{cases} \Omega_{1}^0 = \{1,2,3\}, \ \ \Omega_{2}^0 = \emptyset \\ \Omega_{1}^1 = \{1,3\}, \ \ \Omega_{2}^1 = \{2\} \\ \Omega_{1}^2 = \{1\}, \ \ \Omega_{2}^2 = \{2,3\}   \end{cases} \ \ \text{and} \\[1mm]
		\begin{tabular}{|P{1cm}|P{1cm}|P{1cm}|P{1cm}|}
		\hline 
		$r_{q,\ell}$ & $\ell = 0$ & $\ell = 1$ & $\ell = 2$ \\
		\hline 
		$q = 1$ & 3 & 2 & 1 \\
		\hline
		$q = 2$ & 3 & 3 & 3 \\
		\hline
		$q = 3$ & 3 & 2 & 2 \\
		\hline
		\end{tabular} \hspace{15mm}
		\end{align*}
		The values $r_{q,\ell}$ represent the intermediate reduction orders for each subsystem. Moreover, the transition of the diagonal Gramians $_\ell \hat{\bLa}_q$ for $\ell \in \{0,1,2\}$ is made as follows:
		
		\vspace{2mm}
		
		\noindent
		Step $\ell = 0 \rightarrow$  At this step,  write the original balanced Gramians $_0 \hat{\bLa}_q = \bLa_q, \ \ q \in \Omega$.
		\begin{align*}
		_0 \hat{\bLa}_1 = \left[ \begin{array}{ccc}
		\sigma_{1,1} & 0 & 0  \\ 0 & \sigma_{1,2} & 0  \\ 0 & 0 & \sigma_{1,3} 
		\end{array} \right], \ _0 \hat{\bLa}_2 = \left[ \begin{array}{ccc}
		\sigma_{2,1} & 0 & 0  \\ 0 & \sigma_{2,2} & 0  \\ 0 & 0 & \sigma_{2,3}
		\end{array} \right], \ _0 \hat{\bLa}_3 = \left[ \begin{array}{ccc}
		\sigma_{3,1} & 0 & 0  \\ 0 & \sigma_{3,2} & 0  \\ 0 & 0 & \sigma_{3,3} 
		\end{array} \right].
		\end{align*}
		\noindent
		Step $\ell = 1 \rightarrow$  \ Error bound: $\Vert \hat{\by}_{0} - \hat{\by}_1 \Vert_{2} \leqslant 2 \max(\sigma_{1,3},\sigma_{3,3}) \Vert \bu \Vert_{2}$. \\
		\begin{align*}
		_1 \hat{\bLa}_1 = \left[ \begin{array}{cc}
		\sigma_{1,1} & 0  \\ 0 & \sigma_{1,2} 
		\end{array} \right], \ _1 \hat{\bLa}_2 = \left[ \begin{array}{ccc}
		\sigma_{2,1} & 0 & 0  \\ 0 & \sigma_{2,2} & 0  \\ 0 & 0 & \sigma_{2,3}  
		\end{array} \right], \ _1 \hat{\bLa}_3 = \left[ \begin{array}{cc}
		\sigma_{3,1} & 0  \\ 0 & \sigma_{3,2}   
		\end{array} \right].
		\end{align*}
		Step $\ell = 2\rightarrow$ \ Error bound: $\Vert \hat{\by}_{1} - \hat{\by}_2 \Vert_{2} \leqslant 2 \sigma_{1,2} \Vert \bu \Vert_{2}$. \\
		\begin{align*}
		_2 \hat{\bLa}_1 = 
		\sigma_{1,1} , \ _2 \hat{\bLa}_2 = \left[ \begin{array}{ccc}
		\sigma_{2,1} & 0 & 0  \\ 0 & \sigma_{2,2} & 0  \\ 0 & 0 & \sigma_{2,3}  
		\end{array} \right], \ _2 \hat{\bLa}_3 = \left[ \begin{array}{cc}
		\sigma_{3,1} & 0  \\ 0 & \sigma_{3,2}   
		\end{array} \right].
		\end{align*}
		By combining the two inequalities from steps 1 and 2, it follows that
		\begin{equation*}
		\Vert \by - \hat{\by} \Vert_{2}  \leqslant 2 \big{(} \max(\sigma_{1,3},\sigma_{3,3}) + \sigma_{1,2} \big{)} \Vert \bu \Vert_{2}.
		\end{equation*}
	}
\end{example}

\normalsize

\subsubsection{Stability preservation}

Stability preservation is a very sought after property when devising MOR techniques. As pointed out in \cite{book1LSS}, a switched system is stable if all individual subsystems are stable and the switching is sufficiently slow to permit the transient effects to vanish after each switching time. In this book, Chapter 3.2 presents stability under slow switching with multiple Lyapunov functions.

We present a definition of stability in a uniformly exponentially sense and with imposing again the condition of a minimal dwell time $\mu$. This definition was initially introduced in \cite{book1LSS}.
Moreover, we will show that the reduced order models constructed through the proposed balancing reduction technique, satisfy the conditions of this particular type of stability.

\begin{definition}
	A linear switched system $\Si$ as described in (\ref{LSS_def}),  is uniformly exponentially stable with dwell time $\mu$ if there exist constants $K,M >0$ such that for any solution $(\bx,\bu,\sigma,\by)$, the inequality holds for any $t \geqslant 0$,
	\begin{equation}\label{def_exp_stab}
	\Vert \bx(t) \Vert_2 \leqslant K e^{-Mt} \Vert \bx(0) \Vert_2.
	\end{equation}
	for a control input considered to be zero (i.e. $\bu = \bfz$) and the switching signal $\sigma(t)$ having minimum dwell time $\mu >0$.
\end{definition}

\begin{assumption}
	{\rm
		There exist positive definite matrices $\cQ_i > \bfz, \ i = 1,\ldots,D$ such that
		\begin{equation}\label{lyap_ineq_stabQ1}
		\bA_q^T \cQ_q + \cQ_q \bA_q + 2 M \cQ_q < \bfz, \ \ \forall q \in \Omega
		\end{equation}
		Additionally, assume we can always find positive constants $M, \mu >0$ so that the following inequalities hold for any $q_1,q_2 \in \Omega$,
		\begin{equation}\label{lyap_ineq_stabQ2}
		e^{-M \mu} \bK_{q_1,q_2}^T \cQ_{q_2} \bK_{q_1,q_2} < \cQ_{q_1}.
		\end{equation} }
\end{assumption}

\begin{lemma}
	Consider an LSS $\Si$ for which the conditions in Assumption 45 are satisfied. Then it follows that the system $\Si$ is uniformly exponentially stable with dwell time $\mu$.
\end{lemma}
{\bf{Proof of Lemma 4}}. Let $(\bx,\bu,\sigma,\by)$ be a solution of the LSS with $\bu = \bfz$ and switching signal $\sigma = (q_1,t_1) (q_2,t_2) \ldots$ with minimum dwell time $\mu >0$ (i.e. $t_i \geqslant \mu, \ \forall i$). Again, set $\bV(\bx(t)) = \bx^T(t) \cQ_{q_i} \bx(t), \ \forall t \in (T_{i-1},T_i]$. From (\ref{lyap_ineq_stabQ1}), it directly follows that $\frac{\partial V(\bx(t))}{\partial t} \leqslant -2M V(\bx(t))$. Next, introduce the function 
$$
W(\bx(t)) = e^{2M (t-T_{i-1})} V(\bx(t))  = e^{2M(t-T_{i-1})} \bx^T(t) \cQ_{q_i} \bx(t), \ \forall t \in (T_{i-1},T_i],
$$
and hence, the inequality $\frac{\partial W(\bx(t))}{\partial t} \leqslant 0$ holds. Using the same notations as in (\ref{notationQ}), we get that $W(\bx(t)) \leqslant W(\bx(T_{i-1}^+)) \Rightarrow e^{2M (t-T_{i-1})} V(\bx(t)) \leqslant  V(\bx(T_{i-1}^+))$. Then
\begin{equation}\label{ineqVstab}
V(\bx(t)) \leqslant e^{-2M (t-T_{i-1})}  V(\bx(T_{i-1}^+)), \  \forall t \in (T_{i-1},T_i].
\end{equation}
Now using that $\bx(T_{i-1}^+) = \bK_{q_{i-1},q_{i}} \bx(T_{i-1})$, write
\begin{equation}\label{VTiplusstab}
V(\bx(T_{i-1}^+)) =  \bx^T(T_{i-1}) \bK_{q_{i-1},q_{i}}^T \cQ_{q_{i}} \bK_{q_{i-1},q_{i}} \bx(T_{i-1}).
\end{equation}
From (\ref{lyap_ineq_stabQ2}) and (\ref{VTiplusstab}), we get that
\begin{equation}\label{ineqVTiplusstab}
V(\bx(T_{i-1}^+)) \leqslant  e^{M \mu} \bx^T(T_{i-1}) \cQ_{i-1} \bx(T_{i-1}) =   e^{M \mu} V(\bx(T_{i-1})).
\end{equation}
By plugging in $t = T_i$ in (\ref{ineqVstab}) and using (\ref{ineqVTiplusstab}), it follows that
\begin{equation}\label{VTistab}
V(\bx(T_i)) \leqslant e^{-2Mt_i+M \mu} V(\bx(T_{i-1})).
\end{equation}
By putting all the relations in (\ref{VTistab}) together ($k \in \{1,2,\ldots,i\}$), write that
\begin{align}\label{VTistab2}
V(\bx(T_i)) & \leqslant e^{-M(2t_i-\mu)} V(\bx(T_{i-1})) \leqslant e^{-M(2t_i+2t_{i-1}-2\mu)} V(\bx(T_{i-2})) \nonumber \\ & \leqslant \ldots \leqslant e^{-M(2 T_i-i \mu)} V(\bx(0)).
\end{align}
Since $t > T_{i-1} = \sum_{k=1}^{i-1} t_k$ and by using the fact that the system has minimum dwell time $\mu$ in each operational mode, i.e. $t_k \geqslant \mu$, it follows that $t > (i-1)\mu$. Furthermore, by putting together (\ref{ineqVstab}), (\ref{ineqVTiplusstab}) and (\ref{VTistab2}), the results hold $\forall t \in (T_{i-1},T_i]$,
\begin{align}
V(\bx(t)) \leqslant e^{-2M (t-T_{i-1})}   e^{M \mu} e^{-M( 2T_{i-1}-(i-1)\mu)} V(\bx(0)) \nonumber \\ =  e^{-M (2t-i \mu)}    V(\bx(0)) \leqslant  e^{-M(t-\mu)}   V(\bx(0)). \label{Vtstab}
\end{align}
Additionally, assume that for all $q \in \Omega$, the following inequality holds for $\epsilon, \phi >0$
\begin{equation}\label{boundQstab}
\epsilon^2 \cQ_q \leqslant \bI_{n_q} \leqslant \phi^2 \cQ_q.
\end{equation}
From (\ref{Vtstab}) and (\ref{boundQstab}), it follows that for all $t \in (T_{i-1},T_i]$
\begin{align*}
\Vert \bx(t) \Vert_2^2 &= \bx(t)^T \bx(t) \leqslant \phi^2 \bx(t)^T \cQ_{q_i} \bx(t) = \phi^2 V(\bx(t)) \leqslant  \phi^2 e^{-M(t-\mu)}   V(\bx(0)) \\  &=  \phi^2 e^{-M(t-\mu)}  \bx(0)^T \cQ_{q_1} \bx(0) \leqslant  \frac{\phi^2}{\epsilon^2} e^{M\mu} e^{-Mt}  \Vert \bx(0) \Vert_2^2.
\end{align*}
By choosing $K = \frac{\phi^2}{\epsilon^2} e^{M\mu}$, the result of Lemma 4 is proven (from Definition 15). \sq

\begin{corollary}
	Consider an LSS $\Si$ for which the second condition in Assumption 45 is satisfied and, additionally, there exists $M>0$ so that
	\begin{equation}\label{cond_2mu}
	\bA_q^T \cQ_q + \cQ_q \bA_q + M \cQ_q <\bfz, \ \forall q \in \Omega.
	\end{equation}
	It follows that the system $\Si$ is uniformly exponentially stable with dwell time $2\mu$.
\end{corollary}

\begin{corollary}
	Consider an LSS $\Si$ for which the second condition in Assumption 41 is satisfied and, additionally, there exists $M>0$ so that $\bA_q^T \cQ_q + \cQ_q \bA_q + M \cQ_q <\bfz, \ \forall q \in \Omega$.
	Then it follows that the system $\Si$ is uniformly exponentially stable with dwell time $2\mu$ (where $\mu = - \frac{\ln \gamma}{M}$ for $\gamma = \min\limits_{i,j \in \Omega, \ i \neq j}\gamma_{i,j}$).
\end{corollary}

\begin{corollary}
	If the conditions in Assumption 41 or Assumption 42 hold, then the LSS $\Si$ is uniformly exponentially stable with dwell time $2\mu$ (where $\mu$ is constructed as in the proofs of Lemma 1 or, respectively, Lemma 2).
\end{corollary}

\begin{corollary}
	If the conditions in Assumption 41 or Assumption 42 hold for the original LSS model $\Si$, then the same conditions also hold for the reduced-order LSS model $\bar{\Si}$, as introduced in Defininition 4.1, with the same dwell time $\mu$.
\end{corollary}
{\bf{Proof of Corollary 4}}. Consider that the dwell time $\mu$ is constructed in the proof of Lemma 1, i.e.
$e^{-M \mu} \bK_{q_1,q_2}^T \cQ_{q_2} \bK_{q_1,q_2} < \cQ_{q_1} 
$
and also, the Gramians $\cQ_q$ satisfy $\bA_q^T \cQ_q + \cQ_q \bA_q + M \cQ_q <0$. By multiplying the first inequality with $(\bS_{q_1}^{-1})^T$ to the left and with $\bS_{q_1}^{-1}$ to the right, where $\bS_q$ was defined in (\ref{defSq}), we write
\begin{align}
e^{-M \mu} \Big{(} (\bS_{q_1}^{-1})^T  \bK_{q_1,q_2}^T \bS_{q_2}^{T}  \Big{)}  \Big{(}  (\bS_{q_2}^{-1})^T  \cQ_{q_2} \bS_{q_2}^{-1} \Big{)} \Big{(} \bS_{q_2}  \bK_{q_1,q_2} \bS_{q_1}^{-1} \Big{)} < (\bS_{q_1}^{-1})^T \cQ_{q_1} \bS_{q_1}^{-1} \nonumber \\ \Rightarrow e^{-M \mu} \bar{\bK}_{q_1,q_2}^T \bLa_{q_2} \bar{\bK}_{q_1,q_2} < \bLa_{q_1} \label{ineq_bar1} \ \ \ \ \ \ \ \ \ \ \ \ \ \ \ \ \ \ \ \ \ \ \ \ \ \ \ \ \  
\end{align}
Similarly, by multiplying $\bA_q^T \cQ_q + \cQ_q \bA_q + M \cQ_q <0$ with $(\bS_{q}^{-1})^T$ to the left and with $\bS_{q}^{-1}$ to the right, it follows that
\begin{align}
\Big{(} (\bS_{q}^{-1})^T \bA_q^T \bS_q^T \Big{)} \Big{(} (\bS_{q}^{-1})^T \cQ_q \bS_{q}^{-1} \Big{)} + \Big{(} (\bS_{q}^{-1})^T \cQ_q \bS_{q}^{-1} \Big{)} \Big{(} \bS_q \bA_q \bS_q^{-1} \Big{)} \nonumber \\ + M \Big{(} (\bS_{q}^{-1})^T \cQ_q \bS_{q}^{-1} \Big{)} <0  \Rightarrow \bar{\bA}_q^T \bLa_q + \bLa_q \bar{\bA}_q + M \bLa_q < \bfz. \label{ineq_bar2} 
\end{align}
By plugging in the partitioned matrices from (\ref{partition_bal}) into (\ref{ineq_bar1}), it follows
\begin{align}
& e^{-M \mu} \left[ \begin{array}{cc} \bar{\bK}_{q_1,q_2}^{11} & \bar{\bK}_{q_1,q_2}^{12} \\ \bar{\bK}_{q_1,q_2}^{21} & \bar{\bK}_{q_1,q_2}^{22}
\end{array} \right]^T \left[ \begin{array}{cc} \hat{\bLa}_{q_2} & 0 \\ 0 & \check{\bLa}_{q_2}
\end{array} \right] \left[ \begin{array}{cc} \bar{\bK}_{q_1,q_2}^{11} & \bar{\bK}_{q_1,q_2}^{12} \\ \bar{\bK}_{q_1,q_2}^{21} & \bar{\bK}_{q_1,q_2}^{22}
\end{array} \right]  < \left[ \begin{array}{cc} \hat{\bLa}_{q_1} & 0 \\ 0 & \check{\bLa}_{q_1}
\end{array} \right] \nonumber \\ & 
\Rightarrow e^{-M \mu} \Big{(} \big{(} \hat{\bK}_{q_1,q_2} \big{)}^T \hat{\bLa}_{q_2} \hat{\bK}_{q_1,q_2} + \big{(} \bar{\bK}_{q_1,q_2}^{21} \big{)}^T \check{\bLa}_{q_2} \bar{\bK}_{q_1,q_2}^{21} \Big{)} < \hat{\bLa}_{q_1}  \nonumber \\ & \Rightarrow e^{-M \mu} \big{(} \hat{\bK}_{q_1,q_2} \big{)}^T \hat{\bLa}_{q_2} \hat{\bK}_{q_1,q_2}   < \hat{\bLa}_{q_1}, \label{ineq_hat1}
\end{align} 
where $\hat{\bK}_{q_1,q_2} = \bar{\bK}_{q_1,q_2}^{11}$. also, by substituting the partitioned matrices from (\ref{partition_bal}) into (\ref{ineq_bar2}), write
\begin{align}
\left[ \begin{array}{cc} \bar{\bA}_{q}^{11} & \bar{\bA}_{q}^{12} \\ \bar{\bA}_{q}^{21} & \bar{\bA}_{q}^{22}
\end{array} \right]^T \left[ \begin{array}{cc} \hat{\bLa}_{q_2} & 0 \\ 0 & \check{\bLa}_{q_2}
\end{array} \right] + \left[ \begin{array}{cc} \bar{\bA}_{q}^{11} & \bar{\bA}_{q}^{12} \\ \bar{\bA}_{q}^{21} & \bar{\bA}_{q}^{22}
\end{array} \right] \left[ \begin{array}{cc} \hat{\bLa}_{q_2} & 0 \\ 0 & \check{\bLa}_{q_2}
\end{array} \right]   + M \left[ \begin{array}{cc} \hat{\bLa}_{q_1} & 0 \\ 0 & \check{\bLa}_{q_1}
\end{array} \right] < \bfz \nonumber \\ \Rightarrow \hat{\bA}_q^T \hat{\bLa}_q+ \hat{\bLa}_q \hat{\bA}_q + M \hat{\bA}_q < \bfz. \hspace{25mm} \label{ineq_hat2} 
\end{align}
From (\ref{ineq_hat1}) and (\ref{ineq_hat2}), the results of Corollary 4 directly follow.

\begin{corollary}
	By putting together the results of Corollary 3 and Corollary 4 it follows that, if the original LSS $\Si$ is uniformly exponentially stable with dwell time $\mu$, then the reduced LSS $\hat{\Si}$ for which the diagonal Gramians $\hat{\bLa}_q$ are associated, is also uniformly exponentially stable with dwell time $2\mu$.
\end{corollary}

\section{Numerical Examples}

\subsection{Small system with 3 modes}

Consider the case for which $D = 3$. The reachability Gramians $\cP_i, \ i \in \{1,2,3\}$, satisfy the following equations
\begin{align*}
\bA_1 \cP_1 +\cP_1 \bA_1^T+  \bK_{2,1} \cP_2 \bK_{2,1}^T+  \bK_{3,1} \cP_3 \bK_{3,1}^T + \bB_1 \bB_1^T = \bfz, \\
\bA_2 \cP_2 +\cP_2 \bA_2^T+  \bK_{1,2} \cP_1 \bK_{1,2}^T+  \bK_{3,2} \cP_3 \bK_{3,2}^T + \bB_2 \bB_2^T = \bfz, \\
\bA_3 \cP_3 +\cP_3 \bA_3^T+  \bK_{1,3} \cP_1 \bK_{1,3}^T+  \bK_{2,3} \cP_2 \bK_{2,3}^T + \bB_3 \bB_3^T = \bfz.
\end{align*}
which can be compactly written as
\begin{equation}\label{PgramDLyap_ex}
\bA_{\bD} \bP_\bD + \bP_\bD \bA_\bD^T +\bK_{\scriptsize\reflectbox{\bD}_1} \bP_\bD \bK_{\scriptsize\reflectbox{\bD}_1}^T+\bK_{\scriptsize\reflectbox{\bD}_2} \bP_\bD \bK_{\scriptsize\reflectbox{\bD}_2}^T+ \bB_\bD \bB_\bD^T = \bfz,
\end{equation}
where $\bA_\bD,\bB_\bD$ and $\bP_\bD$ are as in (\ref{notationD}) and also
\begin{equation}\label{defKD12}
\bK_{\scriptsize\reflectbox{\bD}_1} 
=\left[ \begin{array}{ccc}
\bfz & \bK_{2,1} & \bfz \\ \bfz & \bfz & \bK_{3,2} \\ \bK_{1,3} & \bfz & \bfz
\end{array} \right], \ \ \ \bK_{\scriptsize\reflectbox{\bD}_2} 
=\left[ \begin{array}{ccc}
\bfz  & \bfz & \bK_{3,1} \\ \bK_{1,2} & \bfz & \bfz \\ \bfz & \bK_{2,3} & \bfz
\end{array} \right].
\end{equation}
\noindent
Similarly, the observability Gramians $\cQ_i, \ i \in \{1,2,3\}$, satisfy the following equations
\begin{align*}
\bA_1^T \cQ_1 +\cQ_1 \bA_1+  \bK_{1,2}^T \cQ_2 \bK_{1,2}+  \bK_{3,1}^T \cQ_3 \bK_{3,1}^T + \bC_1^T \bC_1 = \bfz, \\
\bA_2^T \cQ_2 +\cQ_2 \bA_2+  \bK_{2,1}^T \cQ_1 \bK_{1,2}+  \bK_{2,3}^T \cQ_3 \bK_{2,3} + \bC_2^T \bC_2 = \bfz, \\
\bA_3^T \cQ_3 +\cQ_3 \bA_3+  \bK_{3,1}^T \cQ_1 \bK_{3,1}+  \bK_{3,2}^T \cQ_2 \bK_{3,2} + \bC_3^T \bC_3 = \bfz,
\end{align*}
which can also be compactly written as
\begin{equation}\label{PgramDLyap_ex}
\bA_{\bD}^T \bQ_\bD + \bQ_\bD \bA_\bD +\bK_{\scriptsize\reflectbox{\bD}_1}^T \bQ_\bD \bK_{\scriptsize\reflectbox{\bD}_1}+\bK_{\scriptsize\reflectbox{\bD}_2}^T \bQ_\bD \bK_{\scriptsize\reflectbox{\bD}_2}+ \bC_\bD^T \bC_\bD = \bfz,
\end{equation}
where $\bA_\bD$, $\bC_\bD$ and $\bQ_\bD$ are block diagonal as in (\ref{notationDgen})  and  $\bK_{\scriptsize\reflectbox{\bD}_i}$ as in (\ref{defKD12}) for $i \in \{1,2\}$. Choose the following system matrices for $\Si$, as
$$
\bA_1 = \left[\begin{array}{ccc} -1 & 0 & 0\\ 0 & -8 & 0\\ 0 & 0 & -5 \end{array}\right], \ \ \bA_2 = \left[\begin{array}{ccc} -2 & 0 & 0\\ 0 & -9 & 0\\ 0 & 0 & -6 \end{array}\right], \ \ \bA_3 = \left[\begin{array}{ccc} -4 & 0 & 0\\ 0 & -3 & 0\\ 0 & 0 & -7 \end{array}\right],
$$
$$
\bB_1 = \left[\begin{array}{c}  1\\ 2\\ -1 \end{array}\right], \ \ \bB_2 = \left[\begin{array}{c} 1\\ -1\\ \frac{3}{2} \end{array}\right], \ \bB_3 = \left[\begin{array}{c} - \frac{1}{2}\\ -2\\ 1 \end{array}\right], \ \begin{cases} \bC_1 = \left[ \begin{array}{ccc} -1 & 1 & \frac{5}{2} \end{array}\right],  \\[2mm] \bC_2 =  \left[\begin{array}{ccc} 1 & 2 & - \frac{7}{2} \end{array}\right], \\[2mm] \bC_3 = \left[\begin{array}{ccc} - \frac{3}{2} & 1 & - \frac{1}{2} \end{array}\right],
\end{cases}
$$
\small
$$
\bM = \left[\begin{array}{ccc} 1 & -1 & 0\\ 0 & 2 & -3\\ 1 & 0 & \frac{1}{2} \end{array}\right], \ \bN = \left[\begin{array}{ccc} 0 & 2 & - \frac{1}{2}\\ 1 & 1 & -1\\ 0 & 0 & -3 \end{array}\right], \ \ \begin{cases} \bK_{1,2} =  \bM/7, \ \ \bK_{2,3} =  \bM/4, \ \ \bK_{3,1} =  \bM/6, \\[2mm] \bK_{2,1} =  \bN/5, \ \ \bK_{3,2} =  \bN/3, \ \ \bK_{1,3} =  \bN/2. \end{cases}
$$
\normalsize
Next, compute the balanced diagonal Gramians $\bLa_i$ as,
\small
$$
\bLa_1 = \left[\begin{array}{ccc} 0.6174 & 0 & 0\\ 0 & 0.0816 & 0\\ 0 & 0 & 0.0419 \end{array}\right],  \bLa_2 = \left[\begin{array}{ccc} 0.4183 & 0 & 0\\ 0 & 0.1514 & 0\\ 0 & 0 & 0.0138 \end{array}\right],   \bLa_3 =  \left[\begin{array}{ccc} 0.3311 & 0 & 0\\ 0 & 0.0948 & 0\\ 0 & 0 & 0.0172 \end{array}\right].
$$
\normalsize
As for Example 4.1, consider the values of the reduced orders for the three subsystems, as $r_1 = 1, \ r_2  = 3$ and $r_3 = 2$. We recover the system matrices of the reduced LSS $\hat{\Si}$ as,
\small
$$
\hat{\bA}_1 = -1.4152, \ \ \hat{\bA}_2 = \left[\begin{array}{ccc}    -7.7330 &  -2.9578  & -1.4537 \\
1.6867 &  -0.9066  & -0.5297 \\
-0.5775 &  1.1507 & -8.3605 \end{array}\right], \ \ \hat{\bA}_3 = - \left[ \begin{array}{ccc}    2.9416 &  0.7103 \\  1.0000  & 5.0427 \end{array}\right], ,
$$
$$   
\hat{\bB}_1 =  -1.3006, \     \hat{\bB}_2 = \left[\begin{array}{c}  -2.4972 \\ 0.0221
\\ -0.0636 \end{array}\right], \  \hat{\bB}_3 = \left[\begin{array}{c}    1.2816
\\  0.2190 \end{array}\right], \ \ \hat{\bC}_1 = 1.2875,  \ \  \hat{\bC}_2 = \left[\begin{array}{c} 2.4992  \\  0.3182  \\  0.2538  \end{array}\right]^T, 
$$
$$
\hat{\bC}_3 = \left[\begin{array}{cc} -1.2857 & -0.5313  \end{array}\right], \ \hat{\bK}_{2,3} = \left[\begin{array}{ccc}   -0.6887  & -0.5866 &  -0.1771 \\
-0.2778  & -0.5806  & -0.0555 \end{array}\right],  \    \hat{\bK}_{3,1} = \left[\begin{array}{cc} -0.3449  &  0.1360 \end{array}\right].
$$
\normalsize

From Example 4.1, it follows that the following bound holds, i.e. $\Vert \by - \hat{\by} \Vert_{2}  \leqslant 2 \big{(} \max(\sigma_{1,3},\sigma_{3,3}) + \sigma_{1,2} \big{)} \Vert \bu \Vert_{2} = 2 ( 0.0816 + 0.0419) = 0.2471 \Vert \bu \Vert_{2} $. 

Consider the switching signal $\sigma(t)$ depicted in Fig.\;1, which is characterized by the sequence of elements $(1,t_1)(3,t_2)(1,t_3)(2,t_4) \ldots (2,t_9)(3,t_{10})$ with dwell times $t_0 = 0s$ and $t_{10} = 15s$.

By choosing the control input as $\bu(t) =  1/2 \sin(20t)e^{-t/2}+1/20 e^{-t/2}$, and performing a time domain simulation, we display in Fig.\;1, the outputs of the original and reduced systems $\Si$ and $\hat{\Si}$.

\begin{figure}[h] \label{fig1} \vspace{-1mm}
	\begin{center}
		\includegraphics[scale=0.26]{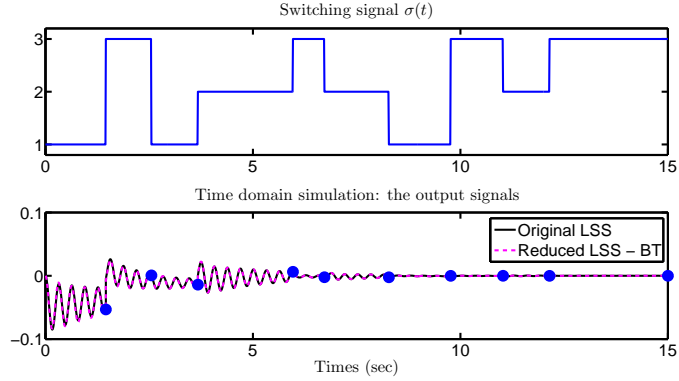}
		\vspace{-3mm}
		\caption{Switching signal $\sigma(t)$ and output $\by(t)$ corresponding to both $\Si$ and $\hat{\Si}$}
	\end{center} \vspace{-3mm}
\end{figure}

\noindent
The absolute value of the difference between the two outputs is presented in Fig.\;2.

\begin{figure}[h] \label{fig2} \vspace{-1mm}
	\begin{center}
		\includegraphics[scale=0.26]{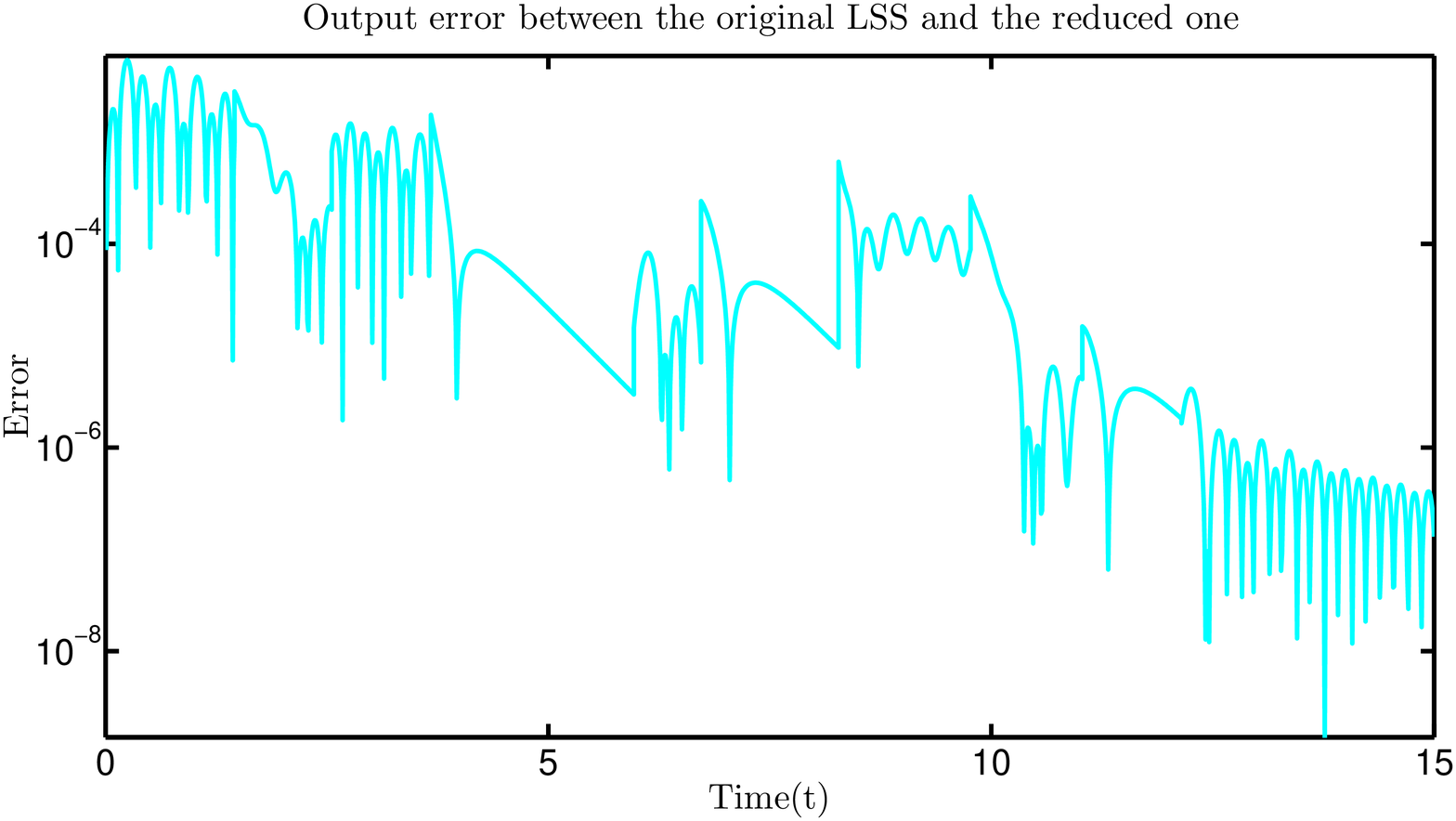}
		\vspace{-3mm}
		\caption{Absolute value of the output error: $\vert \by-\hat{\by} \vert$ }
	\end{center} \vspace{-3mm}
\end{figure}

\subsection{Second example}

For the next experiment, consider the CD player system from the SLICOT benchmark examples for MOR (see \cite{cvd02}). This linear system of order 120 has two inputs and two outputs. We consider that, at any given instance of time, only one input and one output are active (the others are not functional due to mechanical failure). For instance, consider mode $j$ to be activated whenever the $j^{th}$ input and the $j^{th}$ output are simultaneously failing (where $j \in \{1,2\})$.

In this way, we construct an LSS system with two operational modes. Both subsystems are stable SISO linear systems of order 120. This initial linear switched system $\Si$ will be 
reduced by means of the new balanced truncation procedure to obtain $\hat{\Si}_{BT_1}$ and  also by means of the balancing  method proposed in \cite{mtc12} to obtain $\hat{\Si}_{BT_2}$.

There, it has been shown that, if certain conditions
are satisfied, a simultaneous balanced truncation technique can be applied to LSS.  In most practical examples, the existence of a global transformation matrix is not guaranteed. Hence, in \cite{mtc12}, the authors propose instead a method of balancing the so-called average Gramians, i.e. $\cP_{avg} = \frac{1}{D} \sum_{i=1}^D \cP_i$ and $\cQ_{avg} = \frac{1}{D} \sum_{i=1}^D \cQ_i$.

The frequency response of each original subsystem is depicted in Fig.\;3.

\begin{figure}[h] \label{fig3} \vspace{-3mm}
	\begin{center}
		\includegraphics[scale=0.26]{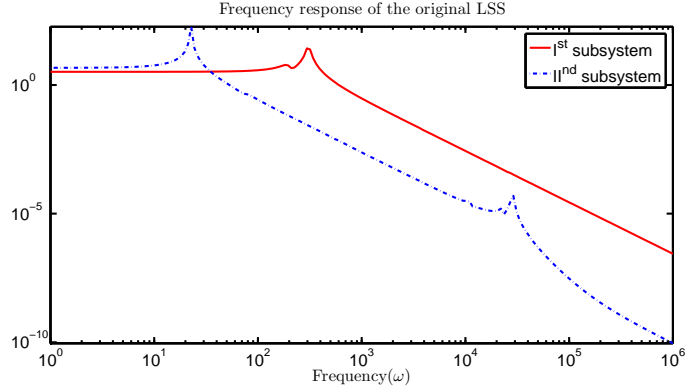}
		\vspace{-3mm}
		\caption{Frequency response of the original subsystems}
	\end{center} \vspace{-6mm}
\end{figure}

Choose the truncation orders $k_1 = k_2 = 33$ for the reduced systems using both methods.
As for the first example, compare the time domain response of the original linear switched system against the ones corresponding to the two reduced models. We use he same signal as in Section 5.1 as control input, i.e. $\bu(t) =  1/2 \sin(20t)e^{-t/2}+1/20 e^{-t/2}$. The switching times $t_i$ are randomly chosen within [0,10]s so that $t_i > 0.5 s, \ \forall i$.

The switching signal $\sigma(t)$ is depicted in the upper part of Fig.\;4, while in the lower part of Fig.\;4, the outputs of the tree LSS mentined above are displayed.  

\begin{figure}[h] \label{fig5} \vspace{-3mm}
	\begin{center}
		\includegraphics[scale=0.26]{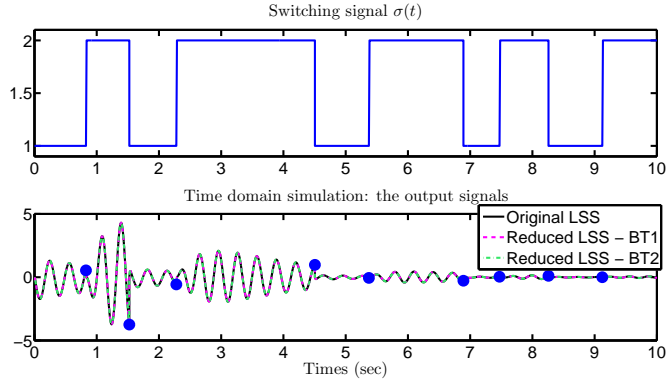}
		\vspace{-3mm}
		\caption{Time domain simulation}
	\end{center} \vspace{-3mm}
\end{figure}

Notice that the output of the original system $\Si$ is well approximated when using any of the two MOR methods.

Finally, by inspecting the time domain error between the original response and the one corresponding to the two reduced models (depicted in Fig.\;5), observe that the new proposed method generally produces better results. The error curve corresponding to the BT1 method is below the error curve corresponding to the BT2 method for most of the points on the time axis.

\begin{figure}[h] \label{fig6} \vspace{-1mm}
	\begin{center}
		\includegraphics[scale=0.26]{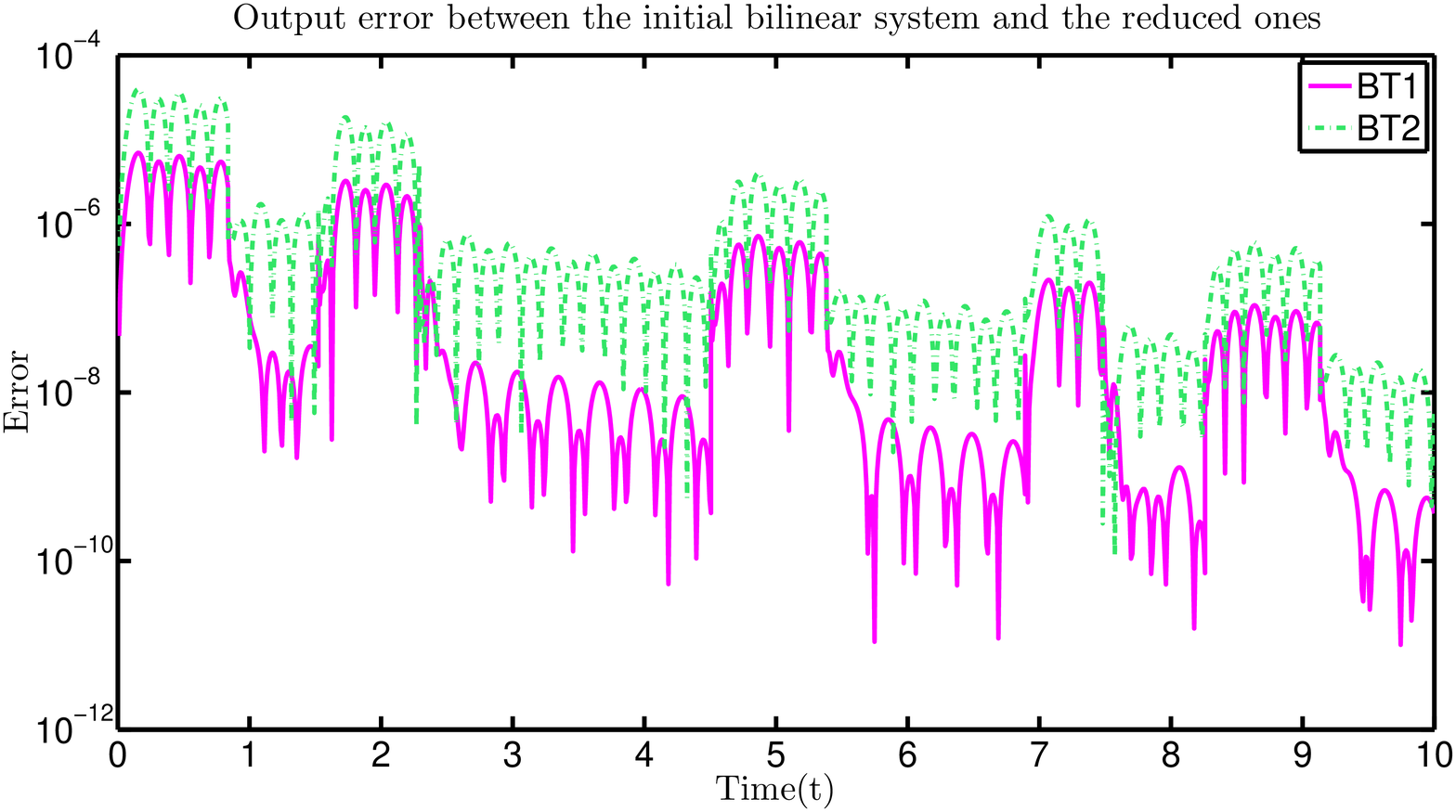}
		\vspace{-3mm}
		\caption{Time domain approximation error}
	\end{center}
\end{figure}

We conclude that the new proposed balancing method produces better results than the one proposed in \cite{mtc12}, in the sense that the original output is better approximated for this particular choice of LSS and control input. Moreover, our method can be applied to LSS with subsystems having different dimensions $n_i, i \in \Omega$ and provide reduced order models again with possibly different dimensions $r_i, i \in \Omega$ in different modes. The other method is constrained to having $n_1 = \ldots = n_D$ so that the computation of the average Gramians $\cP_{av}$ and $\cQ_{av}$ is possible. Also, for BT2 it is assumed that a common Lyapunov function exists, which is arguably restrictive. Moreover, another advantage is that one can derive an error bound of the output error for the new proposed method, as presented in Section 4.2.1. This is also true for the second method proposed in \cite{mtc12}.

\section{Conclusion}

In the current work, we have proposed a balanced truncation procedure for the class of linear switched systems which is based on the computation of infinite energy Gramians. These special matrices can be computed by solving generalized Lyapunov equations instead of solving systems of LMIs. The new balancing method has several advantages.

We provided connections between the new Gramians and system theoretical quantities (observation and controlling energy), by means of lower or upper bounds. Moreover, it turned out that an error bound involving the inputs, outputs and the truncated entries of the Gramians, could be derived. Finally, by applying the proposed procedure, the reduced order LSS can be proven to be uniformly exponentially stable with certain minimum dwell time, given that the original LSS also had this property.

\bibliographystyle{plainurl}

\bibliography{LSS_bal}

\end{document}